\documentclass{manuscript}
\journal{Elsevier}

\begin{document}

\definecolor{mypurple}{rgb}{0.4,0,1}
\definecolor{myred}{rgb}{0.9,0.25,0.36}

\begin{frontmatter}

\title{A fast Chebyshev method for the Bingham closure with application to active nematic suspensions}
\author[cims]{Scott Weady\corref{mycorrespondingauthor}}
\cortext[mycorrespondingauthor]{Corresponding author}
\ead{scott.weady@nyu.edu}
\author[cims,ccb]{Michael J. Shelley}
\author[ccb]{David B. Stein}
\address[cims]{Courant Institute of Mathematical Sciences, New York University, New York, NY, 10012, USA}
\address[ccb]{Center for Computational Biology, Flatiron Institute, Simons Foundation, New York, NY, 10010, USA}

\begin{abstract}
Continuum kinetic theories provide an important tool for the analysis and simulation of particle suspensions. When those particles are anisotropic, the addition of a particle orientation vector to the kinetic description yields a $2d-1$ dimensional theory which becomes intractable to simulate, especially in three dimensions or near states where the particles are highly aligned. Coarse-grained theories that track only moments of the particle distribution functions provide a more efficient simulation framework, but require closure assumptions. For the particular case where the particles are apolar, the Bingham closure has been found to agree well with the underlying kinetic theory; yet the closure is non-trivial to compute, requiring the solution of an often nearly-singular nonlinear equation at every spatial discretization point at every timestep. In this paper, we present a robust, accurate, and efficient numerical scheme for evaluating the Bingham closure, with a controllable error/efficiency tradeoff. To demonstrate the utility of the method, we carry out high-resolution simulations of a coarse-grained continuum model for a suspension of active particles in parameter regimes inaccessible to kinetic theories. Analysis of these simulations reveals that inaccurately computing the closure can act to effectively limit spatial resolution in the coarse-grained fields. Pushing these simulations to the high spatial resolutions enabled by our method reveals a coupling between vorticity and topological defects in the suspension director field, as well as signatures of energy transfer between scales in this active fluid model.
\end{abstract}


\end{frontmatter}

\section{Introduction}

Suspensions of rod-like particles form a broad class of complex fluids. Liquid crystal polymer solutions are one such example, where passive elongated particles, like the tobacco mosaic virus \cite{fraden1989isotropic}, are translated and reoriented by the fluid, modifying its rheological properties \cite{Feng:1999,Sgalari:2002}. In more recent settings however, the suspended particles generate active stresses through propulsive mechanisms \cite{Dombrowski:2004,Sokolov:2007}, chemically induced surface flows \cite{Wang:2015,Davies-Wykes:2016}, or active cross-linking \cite{Gardel:2004,Koenderink:2009,Kohler:2011}. Such active suspensions can exhibit collective flows at scales orders of magnitude larger than those of the constitutive particles. These large-scale flows, sometimes called active turbulence, are characterized by unsteady, roiling states filled with jets, vortices, and topological defects \cite{SS:2013,Thampi:2016,Doostmohammadi:2018,Duclos:2020}. These compelling non-equilibrium structures have motivated various theoretical models ranging from the particle to continuum levels \cite{Simha:2002,H-OSG:2005,SS:2007,SS:2008,Baskaran:2009}.

Particle-based models provide detailed information but are computationally intractable when the number of particles is large. Continuum kinetic theories provide a powerful alternative to discrete models in the large particle number limit. Here the suspension is represented by means of a particle distribution function which evolves through a nonlinear partial differential equation (PDE) -- a Fokker-Planck equation -- allowing the use of well-established analytical tools and numerical methods \cite{SS:2008,ESS:2013}. Its coefficients are usually grounded in modeling of the microscopic physics. Though less demanding than discrete models, kinetic theories are not immune from computational challenges. Typically the distribution function depends on both particle position and orientation as independent variables, meaning there are $2d-1$ degrees of freedom, with $d$ the spatial dimension. This cost can be reduced by coarse-graining, in which the suspension is represented by macroscopic fields derived from the distribution function. However, the equations of motion for the coarse-grained fields depend on unknown fields which must be approximated through a closure model. Constructing an accurate closure model is therefore essential for preserving the multi-scale dynamics and capturing the correct physics \cite{FCL:1998}.

  Closures have long been used for computational models in rheology \cite{Ottinger:2009}, many-particle systems \cite{Levermore:1997}, and classical turbulence \cite{Durbin:2018}, and have more recently been applied to active fluids \cite{Woodhouse:2012,Gao:2017b,Chen:2018,Theillard:2019}. Not only do these models provide efficient computational frameworks, but they can also offer alternative analytical approaches \cite{Han:2015,Li:2015}. This paper is concerned with a specific closure model for apolar suspensions called the \textit{Bingham closure}, originally introduced by Chaubal and Leal in the context of liquid crystal polymers \cite{CL:1998}. In this context, only the zeroth and second moments, with respect to the orientation variables, of the distribution function are evolved in time. The fourth moment, which appears in the corresponding evolution equations, is then approximated as the fourth moment of the Bingham distribution on the unit sphere \cite{Bingham:1974}, whose parameters are computed at each point in space by constraining the zeroth and second moments. (Because the system we consider here is apolar, odd moments do not occur in the dynamics, however other theories may include such moments.) The Bingham closure demonstrates excellent analytical and numerical agreement with the underlying kinetic theory, capturing the same linear instabilities and topological properties of the director field \cite{GBJS:2017}. Accurately computing the closure is essential for maintaining these features, however previous approaches lack robust methods for doing so.

  Existing methods for computing the Bingham closure typically use low order polynomial interpolants from the second to fourth moment tensors, whose coefficients are fit from sample values over the physically feasible domain of the second moment tensor. While these methods are fast, they have limited accuracy. In this paper, we propose a fast Chebyshev method for computing the Bingham closure which maintains the efficiency of polynomial interpolation while achieving near machine precision. The method relies on transforming the domain of eigenvalues of the second moment tensor to a square domain, where the sample points can be chosen on a Chebyshev grid. We numerically compute the Bingham distribution from the second moment tensor on this Chebyshev grid and integrate to obtain the fourth moment tensor. Here we combine spectrally accurate quadrature for the moments with asymptotics to resolve the nearly-singular distribution function at strongly aligned states.
  
  We first restate and coarse-grain a continuum kinetic model for an active suspension and describe how the Bingham closure arises from the coarse-grained theory. We discuss some analytical properties of the closure model, including a proof that it preserves the evolution of the system entropy. We then describe the numerical method for both two- and three-dimensional systems, and evaluate its accuracy and efficiency. This analysis shows that inaccurate computation of the closure reduces the effective spatial resolution, limiting stability and convergence of the underlying numerical method as well as the accessible parameter regimes. Analytical arguments quantify the computational savings of the Bingham closure versus the kinetic theory, which shows impractically high cost for the kinetic theory at strong nematic alignment.  Though we focus on a particular active fluid model, the methods and analyses here equally apply to other apolar kinetic theories, such as those for passive liquid crystal polymer solutions. We conclude with high resolution two- and three-dimensional simulations, focusing on limits of strong steric interactions and large system size. These simulations reveal novel features in the dynamics, including a length scale determination by the steric alignment parameter as well as connections between fluid vorticity and topological defects. 

\section{Model formulation}

  Here we outline a basic model of an active nematic, a more detailed derivation of which can be found in References \cite{ESS:2013} and \cite{GBJS:2017}. Consider a collection of $N$ rod-like particles of length $\ell$ and thickness $b$ such that the aspect ratio is large, ${r = \ell/b \gg 1}$. Each particle generates a surface flow of the form $\U(s) = \sign(s)u_0\p$, where $-\ell/2 \leq s \leq \ell/2$ is the signed arclength along the rod center-line, $u_0$ is the signed surface speed, and $\p$ is the particle's orientation. Because the surface flow is anti-symmetric across $s=0$, the surface flow generates no motion of the particle itself. That is, the particles are immotile and the system is said to be apolar.

  We assume the particles are immersed in a Stokes fluid having linear dimension $L$ and volume $V=L^d$, where $d$ is the spatial dimension. If the number of particles is large, the suspension can be represented by means of a distribution function $\Psi(\x,\p,t)$, which describes the density of particles at center of mass $\x$ with orientation $\p$. Because the number of particles is conserved, this distribution function satisfies a Fokker-Planck equation,
  \begin{equation}
  \frac{\partial\Psi}{\partial t} + \grad\cdot(\dot\x\Psi) + \grad_p\cdot(\dot\p\Psi) = 0,\label{eq:F-P}
  \end{equation}
  where $\grad$ is the spatial gradient and $\grad_p = (\I - \p\p)\cdot\partial/\partial\p$ is the gradient operator projected onto the unit sphere. The conformational fluxes $\dot\x$ and $\dot\p$ in the equation above describe each particle's translational and angular velocities, respectively, and typically depend on the mean field velocity $\u(\x,t)$ and moments of the distribution function. Defining $\langle g(\p) \rangle = \int_{|\p| = 1} g(\p) \Psi ~ d\p$, the relevant moments are the particle concentration $c(\x,t) = \langle 1 \rangle$ and the second-moment tensor $\D(\x,t) = \langle\p\p\rangle$, where $\p\p$ denotes the outer product. The conformational fluxes are then given by
  \begin{align}
  \dot\x &= \u - D_T\grad\log\Psi,\label{eq:xdot}\\
  \dot\p &= (\I - \p\p)\cdot(\grad\u + 2\zeta_0\D)\cdot\p - D_R\grad_p\log\Psi,\label{eq:pdot}
  \end{align}
  where $D_T$ and $D_R$ are the translational and rotational diffusion coefficients. The translational flux (\ref{eq:xdot}) simply says particles move at the local fluid velocity and diffuse, while the rotational flux (\ref{eq:pdot}) represents torques acting on the particles and their rotational diffusion. The torque is generated by the mean-field quantity $\grad\u + 2\zeta_0\D$, which consists of Jeffery's equation modeling particle rotation due to local velocity gradients \cite{Jeffery:1922}, and steric interactions from Maier-Saupe theory, which causes particles to align with the principal axis of $\D$ \cite{MS:1958}. The parameter $\zeta_0$ describes the strength of steric interactions, though its value does not have a precise physical interpretation.

  The biological active fluids we consider typically have small length and velocity scales, so the fluid is well-approximated by the Stokes equation,
  \begin{equation}
  \begin{gathered}
  -\eta\Delta \u + \grad q = \grad\cdot\bm\Sigma,\\
  \grad\cdot\u = 0,
  \end{gathered}\label{eq:stokes}
  \end{equation}
  where $\eta$ is the viscosity, $q(\x,t)$ is the fluid pressure, and $\bm\Sigma(\x,t)$ is the so-called extra stress tensor. The extra stress has three contributions coming from the dipolar active stress due to the surface flow, stress due to particle density and rigidity, and stress caused by steric interactions. Defining the symmetric rate of strain tensor $\E(\x,t) = (\grad\u + \grad\u^T)/2$ and the fourth moment tensor $\S(\x,t) = \langle\p\p\p\p\rangle$, the total stress is given by
  \begin{equation}
  \bm\Sigma = \sigma_a\D + \sigma_c\S:\E - \sigma_s(\D\cdot\D - \S:\D),\label{eq:stress}
  \end{equation}
  where $\sigma_a = -\pi\eta\ell^2 u_0/2\log(2r)$ is the dipole strength, $\sigma_c = \pi\eta\ell^3/6\log(2r)$ arises from particle rigidity, and $\sigma_s = \pi\eta\ell^3\zeta_0/3\log(2r)$ is the strength of steric interactions. Note that the dipole strength $\sigma_a$ has the opposite sign of the imposed surface flow $u_0$. For $u_0 > 0$ the stress is said to be extensile, like that produced by pusher particles, and for $u_0 < 0$ the stress is contractile, like that produced by puller particles. As demonstrated in several studies, the sign of the dipole strength has considerable effects on the system's structure and stability \cite{SS:2008,ESS:2013}.
  
  \subsection{Non-dimensionalization}

  Defining the mean number density $n = N/V$, we choose a reference length $\ell_c = 1/n\ell^2$, velocity scale $|u_0|$, stress scale $\eta|u_0|/\ell_c$, and normalize the distribution function so that $(1/V)\int\int\Psi ~d\p d\x = 1$. In this case the conformational fluxes take the dimensionless form,
	\begin{align}
  \dot\x &= \u - d_T\grad\log\Psi,\label{eq:xdot-nd}\\
  \dot\p &= (\I - \p\p)\cdot(\grad\u + 2\zeta\D)\cdot\p - d_R\grad_p\log\Psi,\label{eq:pdot-nd}
  \end{align}
  and the stress becomes
  \begin{equation}
  \bm\Sigma = \alpha\D + \beta\S:\E - 2\zeta\beta(\D\cdot\D - \S:\D).\label{eq:stress-nd}
  \end{equation}
  Here $\alpha = \sigma_a/\eta|u_0|\ell^2$ is the dimensionless dipole strength, $\zeta = \zeta_0/|u_0|\ell^2$ is the strength of steric interactions, $\beta = \pi n \ell^3/6\log(2r)$ characterizes the density of particles, and $d_T = (n\ell^2/|u_0|)D_T$ and $d_R = (1/n\ell^2|u_0|)D_R$ are the dimensionless translation and rotational diffusion coefficients. The conservation equation (\ref{eq:F-P}) keeps the same form, and the Stokes equation becomes
  \begin{equation}
  \begin{gathered}
  -\Delta\u + \grad q = \grad\cdot\bm\Sigma,\\
  \grad\cdot\u = 0.
  \end{gathered}\label{eq:stokes-nd}
  \end{equation}
  The system of equations (\ref{eq:F-P}) and (\ref{eq:xdot-nd})-(\ref{eq:stokes-nd}) is now a closed system which we call the kinetic theory.

  \subsection{Moment equations}

  The full kinetic theory is complex and high dimensional which makes it expensive to simulate. By taking moments of the Fokker-Planck equation (\ref{eq:F-P}), we can instead represent the dynamics in terms of coarse-grained fields which depend only on space \cite{SS:2013}. Integrating (\ref{eq:F-P}) over the unit sphere in orientation space $\set{\p:|\p|=1}$ leads to an advection-diffusion equation for the particle concentration $c$,
  \begin{equation}
  \frac{\partial c}{\partial t} + \u\cdot\grad c = d_T\Delta c.\label{eq:dc/dt}\end{equation}
  Similarly, multiplying by $\p\p$ and integrating yields an evolution equation for the tensor $\D$,
  \begin{equation}
  \D^\grad + 2\S:\E = 4\zeta(\D\cdot\D - \S:\D)  + d_T\Delta\D - 2dd_R\bpar{\D - \frac{c}{d}\I},\label{eq:dD/dt}
  \end{equation}
	where $\D^\grad = \partial\D/\partial t + \u\cdot\grad\D - (\grad\u\cdot\D + \D\cdot\grad\u^T)$ is the upper-convected time derivative, with the convention $(\grad\u)_{ij} = \partial u_i/\partial x_j$. The second-moment tensor can be used as a quantifier of local alignment. This is more precisely measured by the scalar orientational order parameter,
  \begin{equation}
  s(\x,t) = \frac{d(\mu_1(\x,t) - 1/d)}{d - 1},\label{eq:scalar-order}
  \end{equation}
  where $\mu_1$ is the largest eigenvalue of the normalized second-moment tensor $\Q(\x,t)=\D(\x,t)/c(\x,t)$, often called the tensor orientational order parameter. Notably, the scalar order parameter is zero in an isotropic state $\mu_1 = 1/d$ and unity in a strongly aligned state $\mu_1 = 1$. The eigenvector $\m(\x,t)$ corresponding to the eigenvalue $\mu_1$ is called the director, which can be interpreted as the mean particle orientation. Note that the director is only defined up to sign, and, moreover, is ill-defined in an isotropic state $\D/c = \I/d$.

  Now the PDE (\ref{eq:dD/dt}) depends on the fourth-moment tensor $\S = \langle\p\p\p\p\rangle$ which so far lacks a dynamical equation. One resolution is to take the fourth moment of the Fokker-Planck equation, however this yields an equation that depends on the sixth-moment tensor $\langle\p\p\p\p\p\p\rangle$, posing the same issue. Alternatively, we can approximate $\S$ in terms of the known lower order moments $c$ and $\D$ through a closure model.

  \subsection{The Bingham closure}

  The kinetic theory, being rooted in microscopic modeling, is similar to the classical Doi-Onsager theories for liquid crystal polymers \cite{Doi:1986}. Various closures have been proposed for such theories, a detailed summary of which can be found in \cite{FCL:1998}, however most are not based on self-consistent solutions for the distribution function, but rather on asymptotic or ad-hoc approximations. As a result, such closures fail to reproduce essential properties of the microscopic model.

  Chaubal and Leal introduced a parametric closure scheme for the Doi-Onsager theory which is not only self-consistent, but also yields exact results in the relevant asymptotic regimes \cite{CL:1998}. Specifically, they assume the distribution function takes the form of the Bingham distribution on the unit sphere \cite{Bingham:1974},
  \begin{equation}
  \Psi_B(\x,\p,t) = Z^{-1}(\x,t) e^{\B(\x,t):\p\p},\label{eq:bingham}
  \end{equation}
  where $\B$ is a traceless symmetric tensor and $Z$ is a scalar normalization constant. (Note that, because $\I:\p\p = 1$, translations of the form $\B \mapsto \B + \gamma\I$ only affect the normalization constant $Z$, hence it is sufficient to take $\trace(\B) = 0$ in which case $\B$ is unique \cite{Li:2015}.) The parameters $\B$ and $Z$ can be computed by imposing the moment constraints $c=\langle 1\rangle_B$ and $\D = \langle\p\p\rangle_B$, where $\langle\cdot\rangle_B$ denotes moments of the Bingham distribution. The Bingham closure then consists of an intermediate mapping $\D\mapsto\B[\D]$, after which the distribution function $\Psi_B$ is integrated to obtain the fourth moment tensor,
  \begin{equation} \S_B[\D] = Z^{-1}\int_{|\p|=1} \p\p\p\p ~ e^{\B[\D]:\p\p} ~ d\p.\label{eq:S_B}\end{equation}
  Since $\S_B$ is a function of $\D$, we can re-express the dynamics of our coarse-grained model as a closed system in terms of $c$ and $\D$, with
  \begin{equation}
  \D^\grad + 2\S_B[\D]:\E = 4\zeta(\D\cdot\D - \S_B[\D]:\D)  + d_T\Delta\D - 2dd_R\bpar{\D - \frac{c}{d}\I},\label{eq:dD/dt-bingham}
  \end{equation}
  and the Stokes equation (\ref{eq:stokes-nd}) forced by the extra stress
  \begin{equation}
  \bm\Sigma_B = \alpha\D + \beta\S_B[\D]:\E - 2\zeta\beta(\D\cdot\D - \S_B[\D]:\D).\label{eq:stress-bingham}
  \end{equation}
  
  The Bingham closure has several analytical properties that make it a natural modeling choice. First, the Bingham distribution has a clear physical interpretation, being the unique minimizer of the entropy
  \[ \mathcal S(t) = \int_V\int_{|\p|=1}(\Psi/\Psi_0)\log(\Psi/\Psi_0) ~ d\p d\x\]
  subject to the constraints $c = \langle 1 \rangle$ and $\D = \langle\p\p\rangle$ \cite{Yu:2010}, where $\Psi_0 = 1/2\pi$ (2D) or $1/4\pi$ (3D) is the isotropic distribution function. 
  %
  Moreover, it satisfies the same evolution identity for the system entropy $\mathcal E(t) = \mathcal S(t) + \kappa \mathcal D(t)$,
\begin{equation}
\begin{aligned}
\mathcal E'(t) &= -\frac{d}{\alpha\Psi_0}\bpar{\int_V 2\E:\E  + \beta \E:\S_B:\E ~ d\x}+ (2d\zeta/\Psi_0 - 4dd_R\kappa)\int_V \bpar{\D - \frac{c}{d}\I}:\bpar{\D - \frac{c}{d}\I} ~ d\x\\
&\quad + 8\kappa\zeta\int_V\D:(\D\cdot\D - \S_B:\D)~d\x   -2\kappa d_T \int_V |\grad\D|^2 ~ d\x \\ & \quad -\bbrack{d_T\int_V\int_{|\p|=1}\Psi_B|\grad\log\Psi_B|^2 ~ d\p d\x + d_R\int_V\int_{|\p|=1}\Psi_B|\grad_p\log\Psi_B|^2 ~d\p d\x},
\end{aligned}\label{eq:dS/dt}
\end{equation}
where $\kappa = -d\zeta\beta/2\Psi_0\alpha$ and $\mathcal D(t) = \int (\D-(c/d)\I):(\D-(c/d)\I) ~ d\x$ (see the Appendix for a detailed proof). Importantly, this implies the sources of entropy production and dissipation are equivalent in both models. Further, the Bingham distribution yields exact solutions for both the isotropic and nematic base states, and the linear stability of these base states are in good agreement with the kinetic theory \cite{GBJS:2017}. Finally, coupled with the evolution equation (\ref{eq:dD/dt-bingham}), the Bingham closure preserves the physical trace condition $\trace(\D) = c$, which can be shown by contracting equation (\ref{eq:dD/dt-bingham}) with the identity matrix. Accurately computing the closure is essential for preserving these analytical properties, which is the main objective of the following section.

\section{Numerical method}

  The Bingham closure can be posed as an inverse problem that consists of determining the parameters $\B$ and $Z$ such that the following constraints are satisfied at each point in space,
  \begin{align}
  c(\x,t) &= Z^{-1}(\x,t)\int_{|\p|=1} e^{\B(\x,t):\p\p} ~ d\p,\label{eq:c-constraint}\\
  \D(\x,t) &= Z^{-1}(\x,t)\int_{|\p|=1} \p\p ~ e^{\B(\x,t):\p\p} ~ d\p\label{eq:D-constraint}.
  \end{align}
  As written, this is a $d(d+1)/2 + 1$ nonlinear system for the upper triangular components of the symmetric tensor $\B$ and the normalization constant $Z$. We can use the first equation (\ref{eq:c-constraint}) to solve for $Z$ so that the nonlinear system can be written
  \begin{equation}
  \frac{\D(\x,t)}{c(\x,t)} = \frac{\int_{|\p|=1} \p\p ~ e^{\B(\x,t):\p\p} ~ d\p}{\int_{|\p|=1} e^{\B(\x,t):\p\p} ~ d\p}.\label{eq:Q-constraint}
  \end{equation}
  Further, by rotating into the diagonal frame of $\D$ and using the trace conditions $\trace(\D) = c$ and $\trace(\B) = 0$, this can be reduced to a mapping from the largest $(d-1)$ eigenvalues of $\D/c$ to the largest $(d-1)$ eigenvalues of $\B$, which is a $(d-1)$-dimensional nonlinear system \cite{CL:1998}. For ease of notation, in the following we assume $\D$ and $\S_B$ are normalized by $c$; the argument can be followed identically by replacing $\D\mapsto\D/c$ and $\S_B\mapsto\S_B/c$. 
  
  \subsection{Diagonalization}\label{sec:rotation}

  Because $\D$ is symmetric, it has an eigendecomposition of the form $\D = \bm\Omega\tilde\D\bm\Omega^T$, where $\bm\Omega$ is an orthonormal matrix and $\tilde\D = \text{diag}\set{\mu_i}_{i=1}^d$ is a diagonal matrix consisting of the ordered eigenvalues of $\D$ with $\sum_{i=1}^d \mu_i = 1$. Conjugating the constraint (\ref{eq:D-constraint}) by $\bm\Omega$, we get
  \begin{equation}
  \begin{aligned}
  \tilde\D &= \int_{|\p| = 1} (\bm\Omega^T\p)(\bm\Omega^T\p)^T ~ \Psi_B ~ d\p
  \\ & = \int_{|\p| = 1} \tilde\p\tilde\p^T ~ \tilde \Psi_B ~ d\tilde\p,
  \end{aligned}
  \end{equation}
  where $\tilde\p = \bm\Omega^T\p$ and $\tilde\Psi_B = \tilde Z^{-1}e^{\B:(\bm\Omega\tilde\p\tilde\p^T\bm\Omega^T)}$. Note that because $\bm\Omega$ is orthonormal, the transformation $\p\mapsto \tilde\p$ is simply a re-parameterization of the unit sphere.

  A sufficient condition for the off-diagonal terms in the integral above to be zero is that the matrix $\B$ is also diagonalized by $\bm\Omega$, which is a consequence of the off-diagonal moments being odd in at least one component of the orientation vector $p_i$. (In fact, because $\B$ is unique, this is also a necessary condition.) Under this condition the Bingham distribution takes the diagonal form $\tilde\Psi_B = \tilde Z^{-1} e^{\tilde\B:\tilde\p\tilde\p^T}$, where $\tilde\B = \set{\lambda_i}_{i=1}^d$ is the diagonal matrix of the eigenvalues of $\B$. Using the condition $|\tilde\p| = 1$ and tracelessness of $\B$, we can write
  \begin{equation*}
  \tilde\Psi_B = \tilde Z^{-1}\exp\bpar{\lambda_d + \sum_{k=1}^{d-1}\lambda_k'p_k^2},
  \end{equation*}
  which yields $d-1$ equations for the parameters $\set{\lambda_i' = \lambda_i - \lambda_d}_{i=1}^{d-1}$,
  \begin{equation}
  \mu_i = \frac{\int_{|\p| = 1} p_i^2 \exp\bpar{\sum_{k=1}^{d-1}\lambda_k'p_k^2} ~ d\p}{\int_{|\p| = 1} \exp\bpar{\sum_{k=1}^{d-1}\lambda_k'p_k^2} ~ d\p},\quad i = 1,\ldots,d-1,\label{eq:mu_i}
  \end{equation}
  where we've cancelled the common factor $e^{\lambda_d}$ and dropped tildes in the integrals.

  After solving for the $\lambda_i'$, the rotated fourth-moment tensor $(\tilde\S_B)_{ijk\ell} = \Omega_{mi}\Omega_{nj}\Omega_{pk}\Omega_{q\ell}S_{mnpq}$ can be computed from
  \begin{equation}
  \tilde S_{ijk\ell} = \frac{\int_{|\p| = 1} p_i p_j p_k p_\ell \exp\bpar{\sum_{k=1}^{d-1}\lambda'_kp_k^2} ~ d\p}{\int_{|\p| = 1} \exp\bpar{\sum_{k=1}^{d-1}\lambda'_kp_k^2} ~ d\p},\quad i,j,k,\ell = 1,\ldots,d,\label{eq:ts_iijj}
  \end{equation}
  Because the diagonalized distribution function is even in each $p_i$, only terms of the form $\tilde S_{iijj}$ (which we refer to as the diagonal terms of $\tilde\S_B$) and their permutations are nonzero. We can further simplify computations by taking advantage of the trace identity
  \begin{equation}
  \sums_{k=1}^d \tilde S_{iikk} = \tilde D_{ii} = \mu_i,\label{eq:trace-S}
  \end{equation}
  the last equality of which holds because $\tilde\D$ is diagonalized, so that only $2d-3$ entries of $\tilde\S$ need to be computed. After computing $\tilde\S_B$, we can determine the contractions $\S_B:\D$ and $\S_B:\E$ using the transformation $\bm\Omega$, which we later describe in more detail.

  Solving the nonlinear system (\ref{eq:mu_i}) requires the computation of several integrals on the $(d-1)$-dimensional unit sphere at each point in the domain. Moreover, as $\mu_1\rightarrow1$, which corresponds to the strongly aligned state $\tilde\Psi(\x,\p,t) = c(\x,t)\delta(\p - \hat\x)$, the system becomes ill-conditioned; in fact $\lambda_1\rightarrow\infty$ as $\mu_1\rightarrow1$. As an alternative, we can compute the bounded mapping $\tilde\D\mapsto\tilde\S_B$ in advance and interpolate at each time step. Chaubal and Leal proposed cubic interpolants for these mappings. However, their interpolants only agree to about 3 digits, with less accuracy near the aligned states \cite{CL:1998}. Here we construct Chebyshev interpolants which resolve the mapping to near machine precision while maintaining low computational cost.
  
  \begin{figure}[t!]
  \centering
  \includegraphics{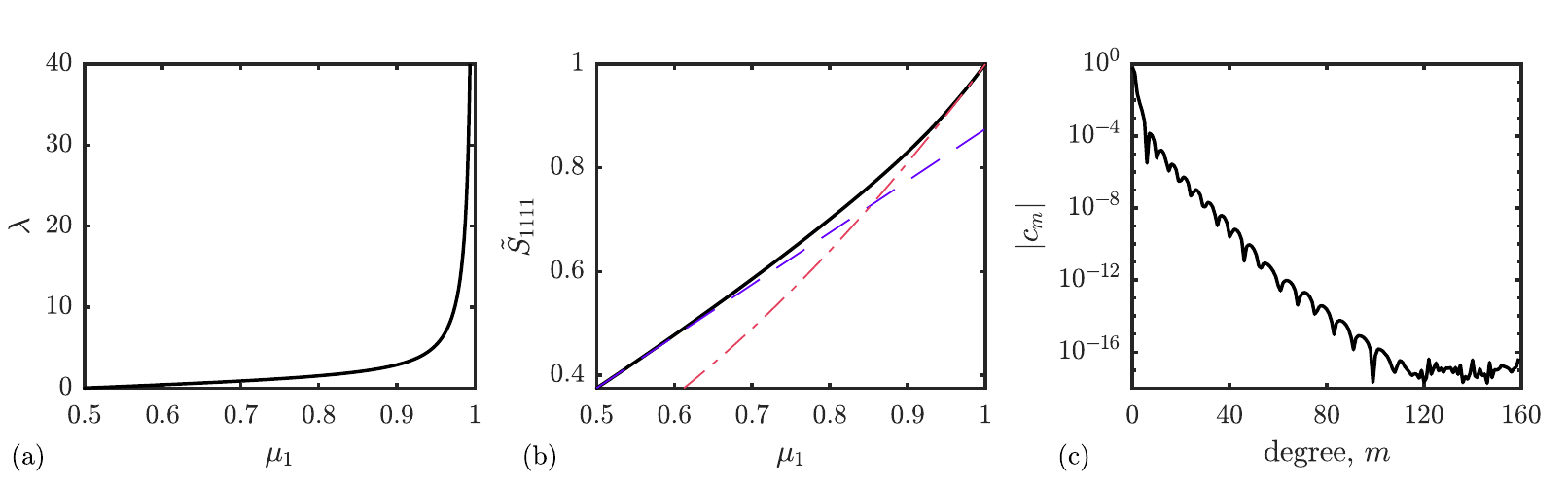}
  \vspace{-0.25in}
  \caption{Two-dimensional Bingham map. Panel (a) shows the parameter $\lambda(\mu_1)$ of the diagonalized Bingham distribution. As $\mu_1\rightarrow 1$, the parameter rapidly increases to infinity which makes inverting equation (\ref{eq:inversion-2d}) ill-conditioned. Panel (b) shows the better-conditioned Bingham map $\mu_1\mapsto\tilde S_{1111}$ with comparison to the linear (\textcolor{mypurple}{\textbf{-- --}}) and quadratic (\textcolor{myred}{\textbf{-- .}}) closures from \cite{SS:2013}. The Bingham map is consistent with each map in their correct limits, matching both pointwise values and first derivatives. Panel (c) shows the corresponding Chebyshev coefficients from Equation \ref{eq:cheb-2d}, where we find approximately 100 modes are needed to resolve the map to near machine precision. (For interpretation of the colors in the figure(s), the reader is referred to the web version of this article.)}\label{fig:bingham-map-2d}
  \end{figure}

  \subsection{Two-dimensional Bingham map}

  Here we detail the construction of the two-dimensional Bingham map. Because $\mu_1 + \mu_2 = 1$, we only need to consider values where the maximum eigenvalue $\mu_1$ is in the interval $[1/2,1]$. In this case the nonlinear system (\ref{eq:c-constraint})-(\ref{eq:D-constraint}) for the largest eigenvalue $\lambda_1 = \lambda$ of the Bingham parameter $\B$ becomes
  \begin{equation}\mu_1 = \frac{\int_0^{2\pi} \cos^2\theta e^{\lambda \cos2\theta} ~ d\theta}{\int_0^{2\pi} e^{\lambda\cos2\theta} ~ d\theta},\label{eq:inversion-2d}\end{equation}
  where we've converted the integrals in Equation (\ref{eq:mu_i}) to polar coordinates with $\p = (\cos\theta,\sin\theta)$. Note that we've used $\lambda_1 = -\lambda_2$ to simplify the exponent.

  The integrals in equation (\ref{eq:inversion-2d}) can be computed analytically using the identity
  \[\int_0^{2\pi}\cos(2n\theta)e^{\lambda\cos2\theta} ~ d\theta = 2\pi I_n(\lambda),\]
  where $I_n(\lambda)$ is the $n$th modified Bessel function of the first kind. We can then write equation (\ref{eq:inversion-2d}) as
  \begin{gather}
   F(\lambda;\mu_1) := \frac{1}{2}\bpar{1 + \frac{I_1(\lambda)}{I_0(\lambda)}} - \mu_1 = 0,\quad\text{with}~ 1/2\leq\mu_1\leq 1.\label{eq:mu1_besseli}
  \end{gather}
  Given $\mu_1 \in [1/2,1]$, we solve this equation for $\lambda(\mu_1)$ with Newton's method, where the Jacobian can similarly be expressed in terms of Bessel functions,
  \[ \frac{\partial F(\lambda;\mu_1)}{\partial \lambda} = \frac{1}{4}\Bigg(1 - 2\Big(\frac{I_1(\lambda)}{I_0(\lambda)}\Big)^2 + \frac{I_2(\lambda)}{I_0(\lambda)}\Bigg).\]
  Once $\lambda(\mu_1)$ has been computed, we evaluate $\tilde S_{1111}(\mu_1)$ using the formula
  \begin{equation}
  \tilde S_{1111}(\mu_1) = \frac{1}{8}\bpar{ 3 + 4\frac{I_1(\lambda(\mu_1))}{I_0(\lambda(\mu_1))} + \frac{I_2(\lambda(\mu_1))}{I_0(\lambda(\mu_1))}}.\label{eq:S1_besseli}
  \end{equation}
  Because the off-diagonal terms of $\tilde\S$ are zero, the physical-frame tensor $\S_B$ can be completely determined by $\tilde S_{1111}$, the trace identity (\ref{eq:trace-S}), and the rotation matrix $\bm\Omega$, so we only need to compute $\mu_1\mapsto\tilde S_{1111}$.
  \subsubsection{Asymptotics near the aligned state $\mu_1\rightarrow 1$}

  Towards the aligned state $\mu_1\rightarrow 1$ we find $\lambda\rightarrow\infty$ for which evaluating the Bessel functions is ill-conditioned. However, only ratios of Bessel functions occur which are bounded for all $\lambda$. To evaluate the ratios when $\lambda \gg 0$, we make use of the series expansion
  \begin{equation} I_n(\lambda) \sim \frac{e^\lambda}{\sqrt{2\pi\lambda}}\sums_{k=0}^\infty \frac{(-1)^k a_k(n)}{\lambda^k},\label{eq:besseli-asym}\end{equation}
  where
  \[ a_k(n) = \frac{\prod_{\ell=0}^k (4n^2 - (2\ell-1)^2)}{8^k k!}.\]
  The ratios can then be stably computed by canceling the leading coefficient $e^\lambda/\sqrt{2\pi\lambda}$ in (\ref{eq:besseli-asym}),
  \begin{equation}
  \frac{I_n(\lambda)}{I_0(\lambda)} \sim \bpar{\sums_{k=0}^\infty \frac{(-1)^k a_k(n)}{\lambda^k}}\bigg/\bpar{\sums_{k=0}^\infty \frac{(-1)^k a_k(0)}{\lambda^k}}.\label{eq:ratio-asym}
  \end{equation}
  We use this asymptotic form whenever $\lambda > 700$, retaining terms up to order $1/\lambda^4$ where the remainder is found to be $O(10^{-15})$. Finally, at the limiting point $\mu_1 = 1$, we set $\tilde S_{1111} = 1$, which is easily shown by taking the limit $\lambda\rightarrow\infty$ in the expressions (\ref{eq:S1_besseli}) and (\ref{eq:ratio-asym}).
  
  \subsubsection{Interpolation}

  We represent the mapping $\mu_1 \mapsto \tilde S_{1111}$ in a Chebyshev basis
  \begin{equation} \tilde S_{1111}(\mu_1) \approx \sums_{m=0}^M c_m T_m(4\mu_1-3),\label{eq:cheb-2d}\end{equation}
  where $T_m(\nu)$ is the $m$th Chebyshev polynomial. To compute the coefficients, we solve equation (\ref{eq:mu1_besseli}) for $\lambda(\mu_1)$ on a Chebyshev grid $\mu_{1,k} = [\cos((2k-1)\pi/2n) + 3]/4$, and evaluate $\tilde S_{1111}(\mu_1) = \tilde S_{1111}(\lambda(\mu_1))$ with equation (\ref{eq:S1_besseli}), using asymptotics when relevant as described above. We then use the MATLAB package \texttt{chebfun} \cite{chebfun} to compute the coefficients $c_m$. To efficiently evaluate this interpolant in practice, we take advantage of the recurrence relation $T_{m+1}(\nu) = 2\nu T_m(\nu) - T_{m-1}(\nu)$ \cite{Broucke:1973}.

  \vspace{0.125in}

  Figure \ref{fig:bingham-map-2d} shows the intermediate map $\mu_1\mapsto \lambda$ and the closure map $\mu_1\mapsto \tilde S_{1111}$ along with the magnitude of its Chebyshev coefficients. We find about 100 modes are needed to resolve the closure map to near machine precision. In Figure \ref{fig:bingham-map-2d}(b), we compare this mapping with two common closures \cite{SS:2013}. These are the linear closure, $\tilde S_{1111} = 3/8 + (\mu_1 - 1/2)$, which is the linear approximation about the isotropic state $\mu_1 = 1/2$, and the quadratic closure $\tilde S_{1111} = \mu_1^2$, which is the correct form in the strongly aligned limit $\mu_1 = 1$. As expected, the Bingham closure produces exact results in both limits, matching not only point-wise values but also derivatives.

  \begin{figure}[t!]
      \centering
      \includegraphics{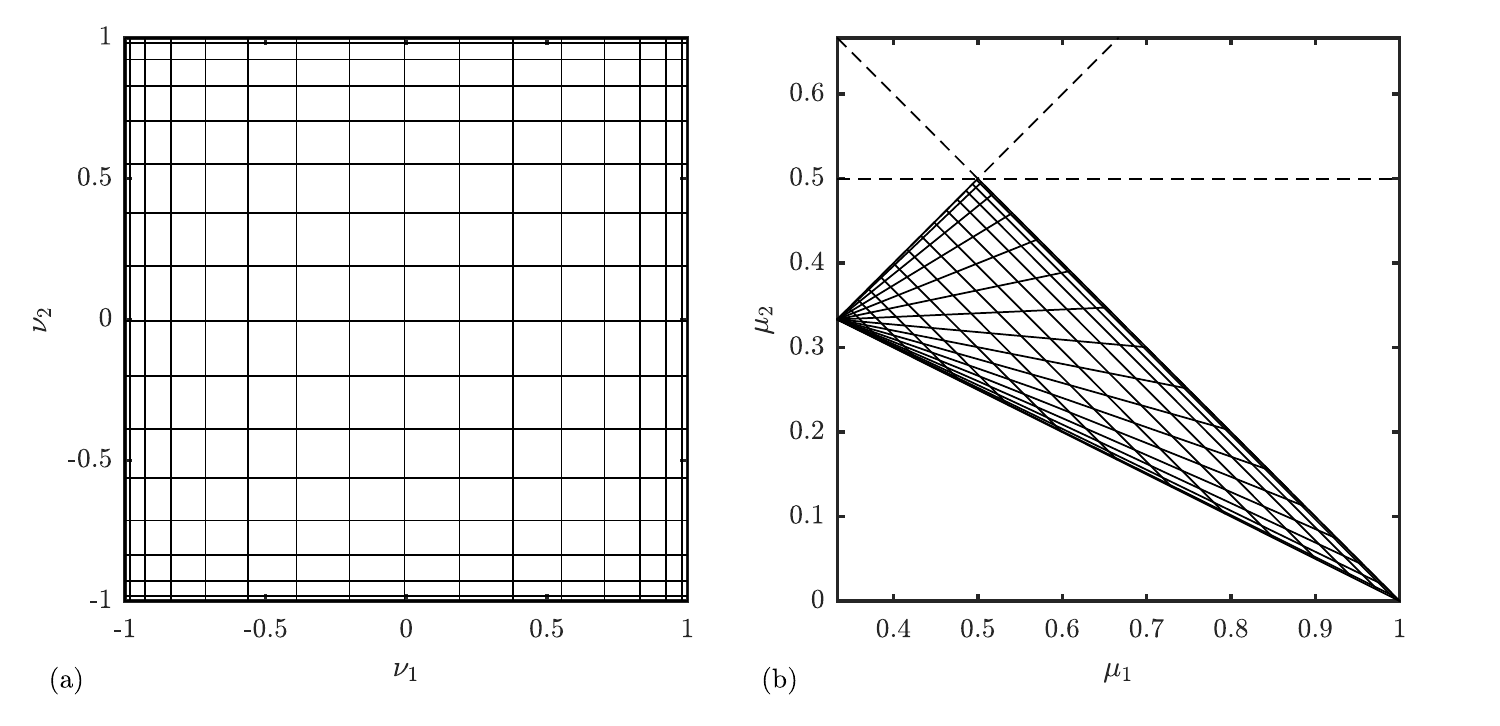}
      \vspace{-0.2in}
      \caption{Mapping $\H^{-1}(\nu_1,\nu_2)$ from the two-dimensional Chebyshev grid $C = \set{(\nu_1,\nu_2):-1\leq\nu_1,\nu_2\leq 1}$ to the feasible domain of eigenvalues $T = \set{(\mu_1,\mu_2):\mu_1+\mu_2+\mu_3 = 1~\text{and}~0\leq\mu_3\leq\mu_2\leq\mu_1\leq1}$. Panel (b) shows the image of the vertical and horizontal contours on the Chebyshev grid, which are well distributed with a tendency to cluster near the isotropic ($\mu_i = 1/3$), planar aligned ($\mu_1 + \mu_2 = 1$), and strongly aligned ($\mu_1 = 1$) states.}
      \label{fig:grid}
  \end{figure}
  
  \subsection{Three-dimensional Bingham map}

  In three dimensions, the physically meaningful domain of the pair of two largest eigenvalues $(\mu_1,\mu_2)$ is determined by the constraints $\mu_1+\mu_2+\mu_3=1$ and $0\leq\mu_3\leq\mu_2\leq\mu_1\leq 1$, which forms a triangle $T$ with corners $(1/3,1/3)$, $(1/2,1/2)$ and $(1,0)$. These corners correspond to the fully isotropic state $(\mu_i = 1/3)$, the planar isotropic state $(\mu_3 = 0)$, and the perfectly aligned state $(\mu_1 = 1)$. From the constraints the boundaries of $T$ are $\mu_1 = \mu_2$, $\mu_2 = \mu_3$, and $\mu_1 + \mu_2 = 1$, the last of which reflects alignment within a plane.
  
  In order to construct a Chebyshev interpolant, we must map $T$ to the square domain $C = \set{(\nu_1,\nu_2):-1\leq\nu_1,\nu_2\leq 1}$. One such mapping can be constructed by composing the linear transformation
  \[ \A(\mu_1,\mu_2) = \begin{pmatrix} \mu_1 - \mu_2 \\ 2\mu_1 + 4\mu_2 - 2 \end{pmatrix}\]
  with the nonlinear transformation
  \[ \G(\mu_1',\mu_2') = \begin{pmatrix} 2(\mu_1' + \mu_2') - 1 \\ \frac{\mu_1' - \mu_2'}{\mu_1' + \mu_2'} \end{pmatrix}. \]
  The resultant $\H = \G\circ\A$ is an invertible mapping from the triangle $T$ to the square $C$. The inverse of this map is simply given by $\H^{-1} = \A^{-1}\circ \G^{-1}$, where
  \[ \G^{-1}(\nu_1,\nu_2) = \begin{pmatrix} \frac{(1 + \nu_1)(1 + \nu_2)}{4} \\ \frac{(1 + \nu_1)(1 - \nu_2)}{4} \end{pmatrix} \]
  and
  \[ \A^{-1}(\nu_1',\nu_2') = \begin{pmatrix} 2\nu_1'/3 + \nu_2'/6 + 1/3 \\ -\nu_1'/3 + \nu_2'/6 + 1/3 \end{pmatrix}.\]
  The image of a separable Chebyshev grid over $C$ under the map $\H^{-1}$ is shown in Figure \ref{fig:grid}. The points are well-distributed across $T$, with clustering near the isotropic, planar aligned, and strongly aligned states. In terms of this transformation, the nonlinear system we need to solve for each $(\nu_1,\nu_2) \in C$ is
  \begin{equation}
  \begin{aligned}
  F_1(\lambda_1,\lambda_2;\nu_1,\nu_2) = \frac{\int_{|\p|=1} p_1^2 e^{\lambda_1 p_1^2 + \lambda_2 p_2^2} d\p}{\int_{|\p|=1} e^{\lambda_1 p_1^2 + \lambda_2 p_2^2} d\p} - H_1^{-1}(\nu_1,\nu_2) = 0,\\
  F_2(\lambda_1,\lambda_2;\nu_1,\nu_2) = \frac{\int_{|\p|=1} p_2^2 e^{\lambda_1 p_1^2 + \lambda_2 p_2^2} d\p}{\int_{|\p|=1} e^{\lambda_1 p_1^2 + \lambda_2 p_2^2} d\p} - H_2^{-1}(\nu_1,\nu_2) = 0.
  \end{aligned}\label{eq:inversion-3d}
  \end{equation}
  Given $(\nu_1,\nu_2) \in C$, we solve this system with Newton's method. Here the Jacobian is
  \begin{equation} \frac{\partial(F_1,F_2)}{\partial(\lambda_1,\lambda_2)} = \begin{pmatrix} \langle p_1p_1p_1p_1 \rangle_B - \langle p_1p_1 \rangle_B \langle p_1p_1 \rangle_B & \langle p_1p_1p_2p_2 \rangle_B - \langle p_1p_1 \rangle_B\langle p_2p_2\rangle_B \\ \langle p_1p_1p_2p_2 \rangle_B - \langle p_1p_1 \rangle_B\langle p_2p_2\rangle_B & \langle p_2p_2p_2p_2 \rangle_B - \langle p_2p_2 \rangle_B\langle p_2p_2\rangle_B\end{pmatrix},\label{eq:jacobian}\end{equation}
  where as before $\langle g(\p) \rangle_B = \int_{|\p|=1} g(\p) \Psi_B ~ d\p$ denotes moments of the Bingham distribution. Note that the Jacobian contains $\tilde S_{1111} = \langle p_1p_1p_1p_1 \rangle_B$, $\tilde S_{1122} = \langle p_1p_1p_2p_2\rangle_B$, and $\tilde S_{2222} = \langle p_2p_2p_2p_2 \rangle_B$ at the converged value. As before, since the off-diagonal terms of $\tilde\S_B$ are zero, we can completely determine $\S_B$ from these three values, the trace identity (\ref{eq:trace-S}), and the rotation matrix $\bm\Omega$.

  \begin{figure}[t!]\centering
  \includegraphics{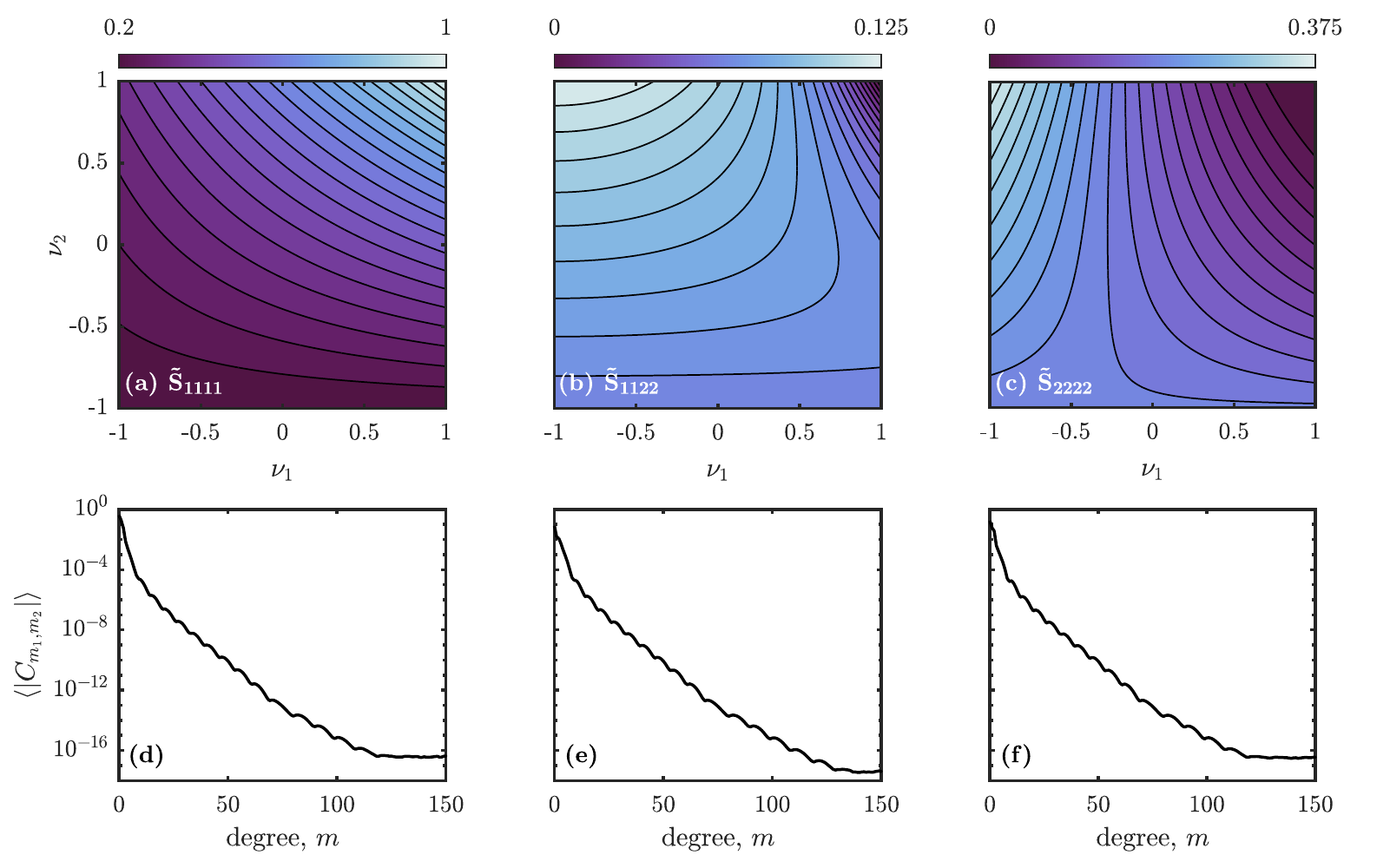}
  \vspace{-0.05in}
  \caption{Three-dimensional Bingham map $(\nu_1,\nu_2)\mapsto(\tilde S_{1111},\tilde S_{1122},\tilde S_{2222})$. Panels (a)-(c) show the computed maps over the transformed domain $(\nu_1,\nu_2) = \H(\mu_1,\mu_2)$. The magnitude of the corresponding Chebyshev coefficients, averaged over each term of degree $m$, are shown below in panels (d)-(f). The coefficients decay near-exponentially, reaching near machine precision at approximately $m = 100$.}\label{fig:bingham-map-3d}
  \end{figure}

  \subsubsection{Quadrature}

  Unlike the two-dimensional case, there is no clear analytical form for the integrals involved in the nonlinear solve. Instead, we compute them numerically in a spherical coordinate system $\p = (\cos\theta,\sin\phi\sin\theta,\cos\phi\sin\theta)$ where $(\phi,\theta) \in [0,2\pi]\times[0,\pi]$, using the spectrally-accurate trapezoidal rule in $\phi$ and Gauss quadrature in $\theta$. Note that we've permuted $p_1$ and $p_3$ from the usual choice of spherical coordinates so that the quadrature nodes cluster at $p_1 = \pm 1$, which is where the Bingham distribution has its peaks in the diagonal coordinate system. For the maps computed here, we used 1024 equispaced nodes in $\phi$ and 4096 Gauss nodes in $\theta$. To avoid overflow for large values of $\lambda_1$ and $\lambda_2$, we subtract $\lambda_1$ from the exponent when numerically evaluating the exponential in the integrand $e^{\lambda_1 p_1^2 + \lambda_2 p_2^2 - \lambda_1}$. Because only ratios of the integrals occur, the common factor $e^{-\lambda_1}$ vanishes and so this does not change the computed moments.

  \subsubsection{Asymptotics at the planar aligned state $\mu_1 + \mu_2 = 1$}

  The nonlinear system (\ref{eq:inversion-3d}) becomes increasingly ill-conditioned as we approach the boundary $\mu_1 + \mu_2 = 1$, or $\mu_3 = 0$. Using the standard spherical coordinate system $\p = (\cos\phi\sin\theta,\sin\phi\sin\theta,\cos\theta)$ and noting that integrand for $\mu_3 = \langle \cos^2\theta \rangle_B$ is strictly positive, at this limit the distribution function must take the form of a $\delta$-function in $\theta$,
  \begin{equation} \Psi_B = Z^{-1}e^{\lambda_1\cos^2\phi\sin^2\theta + \lambda_2\sin^2\phi\sin^2\theta}\delta(\theta-\pi/2).\label{eq:mu3-0}\end{equation}
  Integrating then gives
    \begin{align*} 
    \mu_1 &= \frac{\int_0^{2\pi}\int_0^\pi\cos^2\phi\sin^3\theta e^{\lambda_1'\cos2\phi}\delta(\theta-\pi/2) ~ d\phi d\theta}{\int_0^{2\pi}\int_0^\pi\sin\theta e^{\lambda_1'\cos2\phi}\delta(\theta-\pi/2) ~ d\phi d\theta} 
	 =  \frac{\int_0^{2\pi} \cos^2\phi e^{\lambda_1'\cos2\phi}~d\phi}{\int_0^{2\pi} e^{\lambda'\cos2\phi}~d\phi},
	\end{align*}
	which is the same as the equation for the two-dimensional Bingham map (\ref{eq:inversion-2d}) in the unknown $\lambda_1' = \lambda_1 - \lambda_2$, and can be solved as described before. Because of the constraint $\mu_1 + \mu_2 = 1$, solving for $\lambda_1'$ here is sufficient to determine the full distribution function (\ref{eq:mu3-0}). Once we have $\lambda_1'$, we integrate using the same representation (\ref{eq:mu3-0}) to obtain $\tilde S_{1111},\tilde S_{1122}$, and $\tilde S_{2222}$. In the precomputations, this asymptotic form is used to determine values of the closure map along the boundary $\mu_1+\mu_2 = 1$.

  \subsubsection{Interpolation}

  We represent the maps $(\nu_1,\nu_2)\mapsto\tilde S_{iijj}$ in a separable Chebyshev basis,
 	\begin{equation} \tilde S_{iijj}(\nu_1,\nu_2) \approx \sums_{m_1 + m_2 \leq M} C^{(i,j)}_{m_1m_2}T_{m_1}(\nu_1)T_{m_2}(\nu_2).\label{eq:cheb-3d}\end{equation}
 	As before, $T_m(\nu)$ is the $m$th Chebyshev polynomial, and $\C^{(i,j)}$ is the $(M+1)\times (M+1)$ matrix of coefficients. To compute the expansion, we solve equation (\ref{eq:inversion-3d}) over a two-dimensional Chebyshev grid $(\nu_{1,k},\nu_{2,\ell}) = (\cos((2k-1)\pi/2n),\cos((2\ell-1)\pi/2n))$, after which the coefficients are computed using the extension of \texttt{chebfun} in two-dimensions \cite{chebfun2d}. Note that in practice we must first map the eigenvalues of $\D$ to the transformed domain $\H(\mu_1,\mu_2) = (\nu_1,\nu_2)$, and then evaluate the interpolant in terms of $\nu_1$ and $\nu_2$.

 	\vspace{0.125in}

  Figure \ref{fig:bingham-map-3d} shows the closure maps $(\nu_1,\nu_2)\mapsto(\tilde S_{1111},\tilde S_{1122},\tilde S_{2222})$ over the transformed domain, along with the magnitude of their coefficients, averaged over each term with degree $m$. Similar to the two-dimensional map, about 100 modes are needed to resolve each map to near machine precision. We again make use of the recurrence relation $T_{m+1}(\nu) = 2\nu T_m(\nu) - T_{m-1}(\nu)$ to efficiently evaluate $T_{m_1}(\nu_1)$ and $T_{m_2}(\nu_2)$, each of which only needs to be done once to compute $\tilde S_{1111}(\nu_1,\nu_2),\tilde S_{1122}(\nu_1,\nu_2)$, and $\tilde S_{2222}(\nu_1,\nu_2)$.

  \subsection{Eigendecomposition}\label{section:eigendecomp}

  The diagonal formulation requires an eigendecomposition of the matrix $\D$ at every grid point. Although each matrix is only $d\times d$, computing large numbers of such small decompositions using calls to external routines (e.g. LAPACK) carries significant overhead and complicates efficient parallelization. Here we provide a simple yet robust method to compute the eigendecomposition of $\D$ without external routines in both two and three dimensions.

  \subsubsection{Two dimensions}

  It's straightforward to show the largest eigenvalue $\mu_1$ of $\D$ is given by $2\mu_1 = 1 + \sqrt{2(\D:\D)-1}$, and the corresponding rotation matrix $\bm\Omega$ is
 	\[ \bm\Omega = \begin{pmatrix} \cos\omega & -\sin\omega \\ \sin\omega & \cos\omega \end{pmatrix},\]
  with $2\omega = \arctan[2D_{12}/(2D_{11}-1)]$. Note that the eigenvectors must be arranged as $\bm\Omega = (\v_1 ~ \v_2)$ in descending order.

  \subsubsection{Three dimensions}

  Rather than use explicit formulas for the eigenvalues, which are numerically unstable to evaluate, we instead numerically solve for the roots of the characteristic polynomial of $\D$, which is given by
  \[ p_\D(z) = z^3 - z^2 + a_1 z + a_0,\]
  where $a_1 = -1/2[1 - \trace(\D^2)]$ and $a_0 = -\det(\D)$. (As the eigenvalues are bounded between 0 and 1, this is a well-conditioned problem.) This equation is quickly solved with a few iterations of Newton's method. Once we have one solution $\mu_0$, we can analytically compute the others via
  \begin{equation*}
	\mu_{\pm} = \frac{-(\mu_0-1) \pm \sqrt{(\mu_0-1)^2 - 4(a_1 + \mu_0(\mu_0-1))}}{2}.
	\end{equation*}
  Finally, we sort $\mu_0$ and $\mu_\pm$ so that $\mu_1\geq\mu_2\geq\mu_3$. From the eigenvalues, we can compute the eigenvectors by taking advantage of orthogonality of $\bm\Omega$ \cite{Kopp:2008}. To be concrete, the $i$th eigenvector $\v_i$ with eigenvalue $\mu_i$ satisfies $(\D - \mu_i\I)\cdot\v_i = 0$. Dotting this with the $j$th basis vector $\e_j$ gives $\v_i^T\cdot(\d_j - \mu_i\e_j) = 0$, where $\d_j$ is the $j$th column of $\D$. Notably, this means $\v_i$ is orthogonal to $\d_j - \mu_i\e_j$ for each $j$, which implies
  \begin{equation}\v_i = (\d_1 - \mu_i\e_1)\times(\d_2 - \mu_i\e_2),\quad i = 1,2.\label{eq:eigvec}\end{equation}
  Since the eigenvectors are orthogonal, we can get the final eigenvector by computing the cross product $\v_3 = \v_1\times\v_2$. Computationally, especially near the isotropic state $\D\approx\I/3$, the formula (\ref{eq:eigvec}) yields eigenvectors that are not orthogonal to machine precision, which can result in numerical instability. To stabilize this, we simply redefine $\v_2 = \v_1\times\v_3$ at each point. Finally, we normalize and arrange the eigenvectors in descending order to get the transformation $\bm\Omega = (\v_1 ~ \v_2 ~ \v_3)$\footnote{At the perfectly isotropic state $\D = \I/3$, the method above results in divide by zero errors when normalizing the eigenvectors. We avoid this by perturbing the off-diagonal terms of $\D$ by $10^{-16}$ in all cases.}.

  \subsection{Rotating to the physical frame}\label{sec:inv_rotation}

  The fourth-moment tensor $\S_B$ has many independent components which makes it expensive to store, especially in three dimensions. Fortunately, we only need contractions of $\S_B$ with rank two tensors, that is, $\S_B:\D$ and $\S_B:\E$. This storage can be further reduced by observing that Equations (\ref{eq:stress}) and (\ref{eq:dD/dt}) only depend on the symmetric rank-two tensor $\S_B:\T$ with $\T := \E + 2\zeta\D$. Moreover, we can utilize the rotation-based framework to efficiently compute this contraction. Specifically, write $\S_B$ as the rotation
  \begin{equation} S_{ijk\ell} = \Omega_{im}\Omega_{jn}\Omega_{kp}\Omega_{\ell q} \tilde S_{mnpq},\label{eq:rotation} \end{equation}
  where repeated indices denote summation. We then have
  \begin{align*}
  (\S_B:\T)_{ij} &= (\Omega_{im}\Omega_{jn}\Omega_{kp}\Omega_{\ell q} \tilde S_{mnpq}) T_{k\ell}\\
  & = (\Omega_{im}\Omega_{jn}\tilde S_{mnpq})(\Omega_{kp}\Omega_{\ell q}T_{k\ell})\\
  & = \Omega_{im}\Omega_{jn}(\tilde S_{mnpq} \tilde T_{pq})\\
  & = \Omega_{im}\Omega_{jn}(\tilde\S_B:\tilde\T)_{mn},
  \end{align*}
  or in matrix notation, $\S_B:\T = \bm\Omega(\tilde\S_B:\tilde\T)\bm\Omega^T$ with $\tilde\T = \bm\Omega^T\T\bm\Omega$. Since only the diagonal elements of $\tilde\S_B$ are nonzero, this requires far fewer operations than the explicit rotation formula (\ref{eq:rotation}).

  \subsection{Summary of the closure}

	The method described here has several important features. Most significantly, Chebyshev interpolation preserves the accuracy of a direct nonlinear solve of equation (\ref{eq:mu_i}) at relatively low cost. This interpolation, rather than directly solving the nonlinear system, is essential for efficiency and numerical stability near the aligned state $\mu_1\rightarrow 1$. Further, explicit calculation of the rotations bypasses overhead from eigenvalue routines, which admits efficient parallelization. Lastly, storing and rotating the contraction $\S_B:\T = \bm\Omega( \tilde \S_B:\tilde \T)\bm\Omega^T$ substantially reduces memory requirements and the number of floating point operations. To summarize, the algorithm consists of the following steps:
  \begin{enumerate}[(1)]
    \item At each spatial discretization point, compute the eigendecomposition of the second-moment tensor $\D = \bm\Omega\tilde\D\bm\Omega^T$ using the method described in Section \ref{section:eigendecomp}.
    \item Evaluate the Chebyshev interpolants $\mu_1\mapsto\tilde S_{1111}$ in 2D or $(\mu_1,\mu_2)\mapsto(\tilde S_{1111},\tilde S_{1122},\tilde S_{2222})$ in 3D, and use the trace identities (\ref{eq:trace-S}) to compute the remaining elements of $\tilde \S_B$.
    \item Rotate $\tilde \T = \bm\Omega^T(\E + 2\zeta\D)\bm\Omega$, and compute and store the tensor $\S_B:\T = \bm\Omega(\tilde\S_B:\tilde \T)\bm\Omega^T$.
  \end{enumerate}
  \section{Numerical tests}

  In this section, we evaluate the cost and accuracy of our numerical implementation. We restrict our discussion to the three-dimensional case, finding similar results in two dimensions. The numerical method is based on a pseudo-spectral discretization of Eqns (\ref{eq:stokes-nd}), (\ref{eq:dD/dt-bingham}), and (\ref{eq:stress-bingham}) with the 2/3 anti-aliasing rule, along with a second-order implicit-explicit backward differentiation time-stepping scheme (SBDF2), where the linear terms are handled implicitly and the nonlinear terms explicitly. Both our two- and three-dimensional codes are written in C++ and use OpenMP to parallelize computations. All computations in this section were done on a grid of $256^3$ Fourier modes with a time step $\Delta t = 0.05$.

  For the following tests we set the dimensionless parameters to be $\alpha = -1$, $\beta = 0.8$, and $\zeta = 1$, with box size $L = 15$ and diffusion coefficients $d_T = d_R = 0.1$. This choice of parameters ensures the isotropic and nematic base states are unstable \cite{ESS:2013}, driving persistent chaotic flows, but also guarantees length and time scales are highly resolved for the chosen grid and time-step. In each simulation we initialize $\D$ with a plane-wave perturbation about the isotropic state $\D_0 = \I/3$ such that $\trace(\D) = 1$ and take the concentration to be uniform, $c(\x,t) \equiv 1$. Based on the evolution equation (\ref{eq:dc/dt}), this means the concentration stays uniform for all time.

  \begin{figure}[t!]
  \centering
  \includegraphics{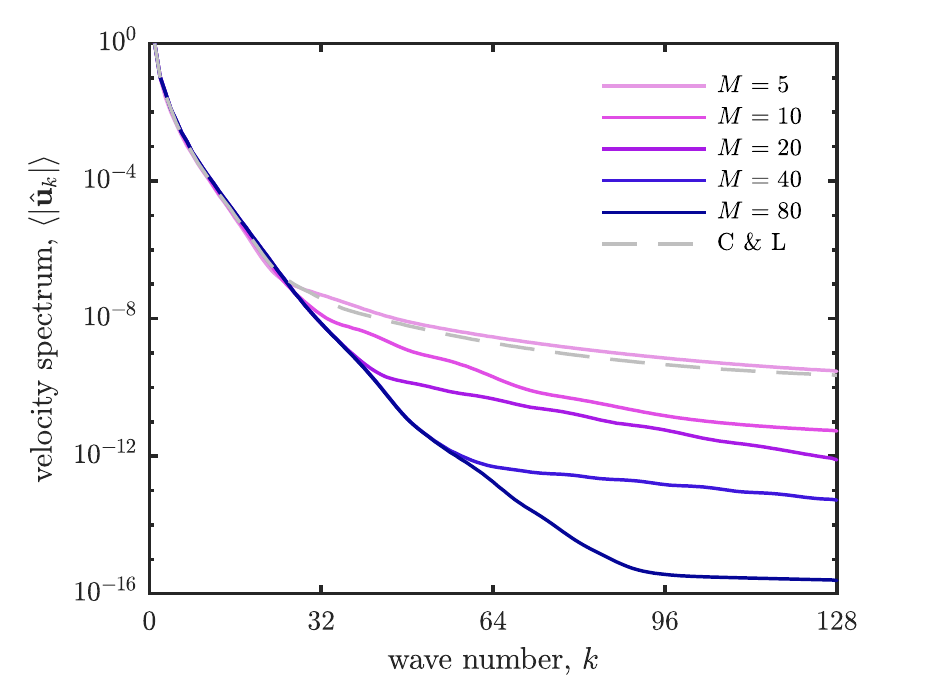}
  \vspace{-0.05in}
  \caption{Convergence in the velocity spectrum $\langle |\hat\u_k| \rangle$, averaged over wave number $|\k| = k$, with the degree $M$ of the Chebyshev interpolant. The gray dashed line indicates the third-order interpolant of Chaubal and Leal (C \& L), for which the velocity is resolved to about $10^{-8}$. This approximately agrees with the fully resolved $M = 80$ interpolant up to $k\approx32$, indicating effectively a fourth of the Fourier expansion is used.}\label{fig:convergence}
  \end{figure}

  \subsection{Spatial convergence}\label{sec:convergence}

  An important and somewhat surprising feature of the Bingham closure is that the accuracy of the entire method is limited by that of the mapping from $\tilde\D$ to $\tilde \S_B$. To demonstrate this, for our initial data and for a fixed value of $M$ of the Chebyshev expansion (\ref{eq:cheb-3d}), we run the simulation to a statistical steady state ($t = 50$). Figure \ref{fig:convergence} shows the resulting velocity spectra $\langle|\hat\u_k|\rangle$ averaged over spherical shells in the wave number $|\k|=k$. As we increase the degree $M$ of the Chebyshev interpolant the dynamical range expands, reaching near machine precision with $M = 80$. For the interpolant of Chaubal and Leal \cite{CL:1998}, the spectrum deviates from the full expansion at wave number $k = 32$ and plateaus, indicating that effectively only a fourth of the potential resolution is used. To be sure this is not an artifact of the eigendecomposition approach in Section \ref{section:eigendecomp}, we performed equivalent calculations in two dimensions using MATLAB's \texttt{eig} function and found the same results.

	The inaccuracy in the velocity field when using a low order interpolant is a consequence of the rotation-based approach in Section \ref{sec:rotation}. To be precise, close to aligned states where the second moment tensor $\D$ has repeated eigenvalues, the eigendecomposition of $\D$ is ill-defined which results in spatial discontinuities in the rotation matrix $\bm\Omega(\x,t)$. In exact arithmetic, these discontinuities are canceled when rotating back to the original frame. However, if the interpolation is not computed accurately, the discontinuities will carry through the inverse rotation. This error becomes even more pronounced in the velocity field due to derivatives of $\S_B$ occurring in the active force $\grad\cdot\bm\Sigma_B$. Note that this problem could be avoided by computing $\S_B$ in the original frame. However, in this case the nonlinear system (\ref{eq:mu_i}) is five-dimensional and defined over an irregular grid, which is not only more expensive and less stable, but also poses further challenges for interpolation.

  \begin{figure}[t!]
  \centering
  \includegraphics{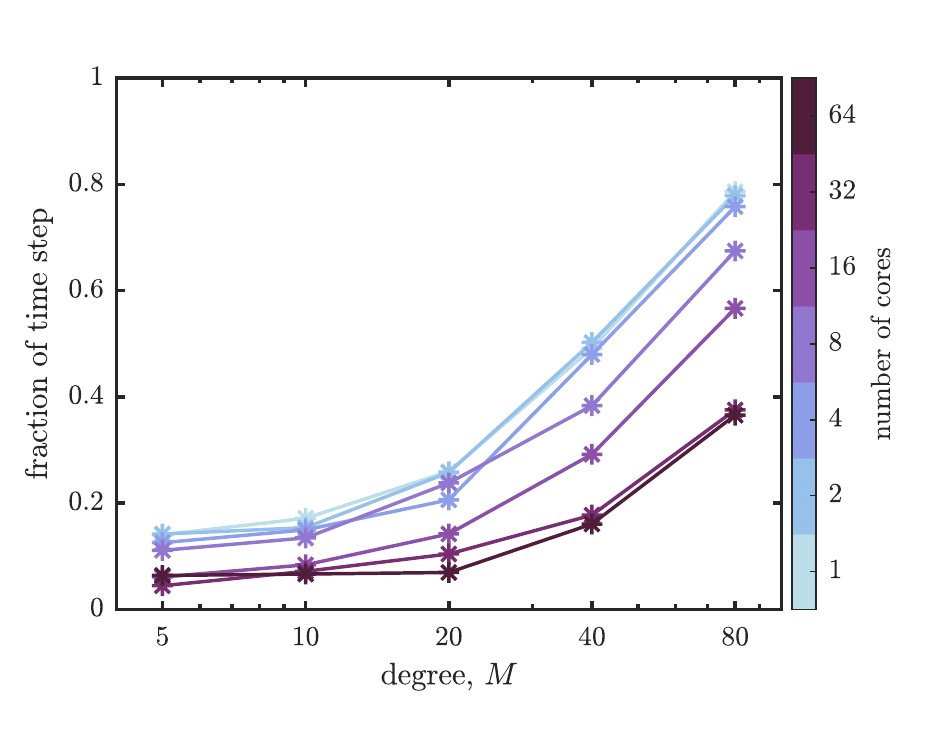}
  \vspace{-0.125in}
  \caption{Fraction of each time step taken by computing the Bingham closure. In serial, the cost between the Bingham closure and the rest of the loop is comparable for $M = 40$. For a fixed degree, the fraction decreases with the number of cores, taking less than half of each time step when the number of cores exceeds 16 for each value of $M$ tested.}\label{fig:timing}
  \end{figure}

  \subsection{Computational cost}

  Representing the closure maps in a Chebyshev series allows us to explicitly balance cost with accuracy. We characterize the cost in Figure \ref{fig:timing}, which shows the fraction of each time step taken by computing the Bingham closure for interpolants of increasing degree. In serial, we find that when $M = 40$ the cost of the Bingham closure is comparable to the remaining cost of each time step. Increasing the number of cores decreases the relative cost, with the Bingham closure taking less than half of a time step for all values of $M$ explored when the number of cores exceeds 16.

  It is also useful to compare the cost of evaluating the Bingham closure through Chebyshev interpolation with that of a direct inversion of equation (\ref{eq:mu_i}). The interpolation has a fixed $O((M+1)^{d-1})$ cost where $M$ is at most 100, while the nonlinear solve is $O(N^{d-1})$, where $N$ is the number of quadrature points in each dimension. We can get an estimate on the number of quadrature points required using the analytical form of the Bingham distribution. Here we make use of identity (\ref{eq:ibp}) from the Appendix,
  \[ \D\cdot\B - \S_B:\B = \frac{d}{2}\bpar{\D - \I/d}.\]
  Conjugating this equation by $\bm\Omega$ gives
  
  \[ \tilde\D\cdot\tilde\B - \tilde\S_B:\tilde\B = \frac{d}{2}\bpar{\tilde\D - \I/d},\]
  which, because $\D,\B$ and $\S_B$ are diagonalized in the same frame, is a diagonal system of equations. Restricting to the two-dimensional case $d = 2$, the first of these equations is
  \[ \mu_1\lambda_1 - (\tilde S_{1111}\lambda_1 + \tilde S_{1122}\lambda_2) = \mu_1 - 1/2. \]
  Using the fact that $\tilde\B$ is trace-free $\lambda_2 = -\lambda_1$ and the trace condition $\tilde S_{1111} + \tilde S_{1122} = \mu_1$, we can solve this equation for $\lambda_1$,
  \[ \lambda_1 = \frac{\mu_1 - 1/2}{2(\mu_1 - \tilde S_{1111})}.\]
  This gives an expression for the standard deviation of the Bingham distribution about $\p = (1,0)$,
  \[ \sigma = \sqrt{\frac{\mu_1 - \tilde S_{1111}}{\mu_1 - 1/2}}.\]
  As an estimate, we demand at least $10$ quadrature nodes within one standard deviation, which, assuming the trapezoidal rule in polar coordinates $\p = (\cos\theta,\sin\theta)$, yields

  \begin{equation} N \approx 40\pi\sqrt{\frac{\mu_1 - 1/2}{\mu_1 - \tilde S_{1111}}}.\label{eq:resolution}\end{equation}

  For example, when $\mu_1 = 0.99$, which regularly occurs when the alignment strength $\xi = 2\zeta/d_R$ is within the physically relevant regime, this estimates $N\approx 345$. In practice, when directly inverting equation (\ref{eq:mu_i}) in a simulation rather than using interpolation, these integrals need to be evaluated several times at every point in space, which may be reasonable in 2D, but is inaccessible with the equivalent estimate in 3D. Moreover, using $\tilde S_{1111} \approx \mu_1^2$ as $\mu_1\rightarrow 1$, we find $N\sim(1 - \mu_1)^{-1/2}$ so that approaching the aligned state requires prohibitive increases in resolution. The estimate (\ref{eq:resolution}) equally applies to the number of discretization points in orientation when simulating the kinetic theory, giving us a rigorous characterization of the savings gained by the closure model. We note that the cost of quadrature can be mitigated by adaptive methods or asymptotic approximations to the moment integrals \cite{Luo:2018}, however such methods are still subject to ill-conditioning of the nonlinear system near the aligned state.

  \section{Numerical simulations}\label{sec:simulations}

    In this section we use the Bingham closure to study two- and three-dimensional suspensions of active extensile particles in the regimes of strong steric interactions and large system size. The Bingham closure is particular useful here as it yields accurate solutions near the isotropic and aligned states which both frequently occur in these regimes. As before, the discretization is pseudo-spectral and we use the implicit-explicit SBDF2 time-stepping scheme, where we use a $4096^2$ grid in two dimensions and a $512^3$ grid in three dimensions, with degree $M = 80$ Chebyshev interpolants in all cases.
    
    \subsection{Strongly aligned dynamics}\label{sec:strongly-aligned}

    Active nematic suspensions exhibit rich topological structures that are an intrinsic part of the system's dynamics \cite{Giomi:2013,Thampi:2014,GBGBS:2015,GBJS:2017,Doostmohammadi:2018}. The primary features are called disclinations, or defects, which refer to points of low orientational order (i.e. the scalar order parameter $s$, defined in equation (\ref{eq:scalar-order}), is approximately zero) at which the director field is ill-defined. In two dimensions, the characteristic topological features are $\pm1/2$ defects, which correspond to a clockwise/counterclockwise rotation of the director about a point of isotropy $s = 0$, respectively. Simulations of a phenomenological Landau-deGennes $Q$-tensor theory have found equivalent features in three dimensions that are closed disclination lines and rings along which the director undergoes various types of three-dimensional rotations \cite{Duclos:2020}. Here we find and examine these topological features using the Bingham closure for the case of strong alignment $\zeta \gg 1$ in both two and three dimensions. The remaining dimensionless parameters are fixed at $\alpha = -1$, $\beta = 0.8$, $d_T = d_R = 0.05$, and $L = 30$. 

\begin{figure}[t!]
  \centering
  \includegraphics{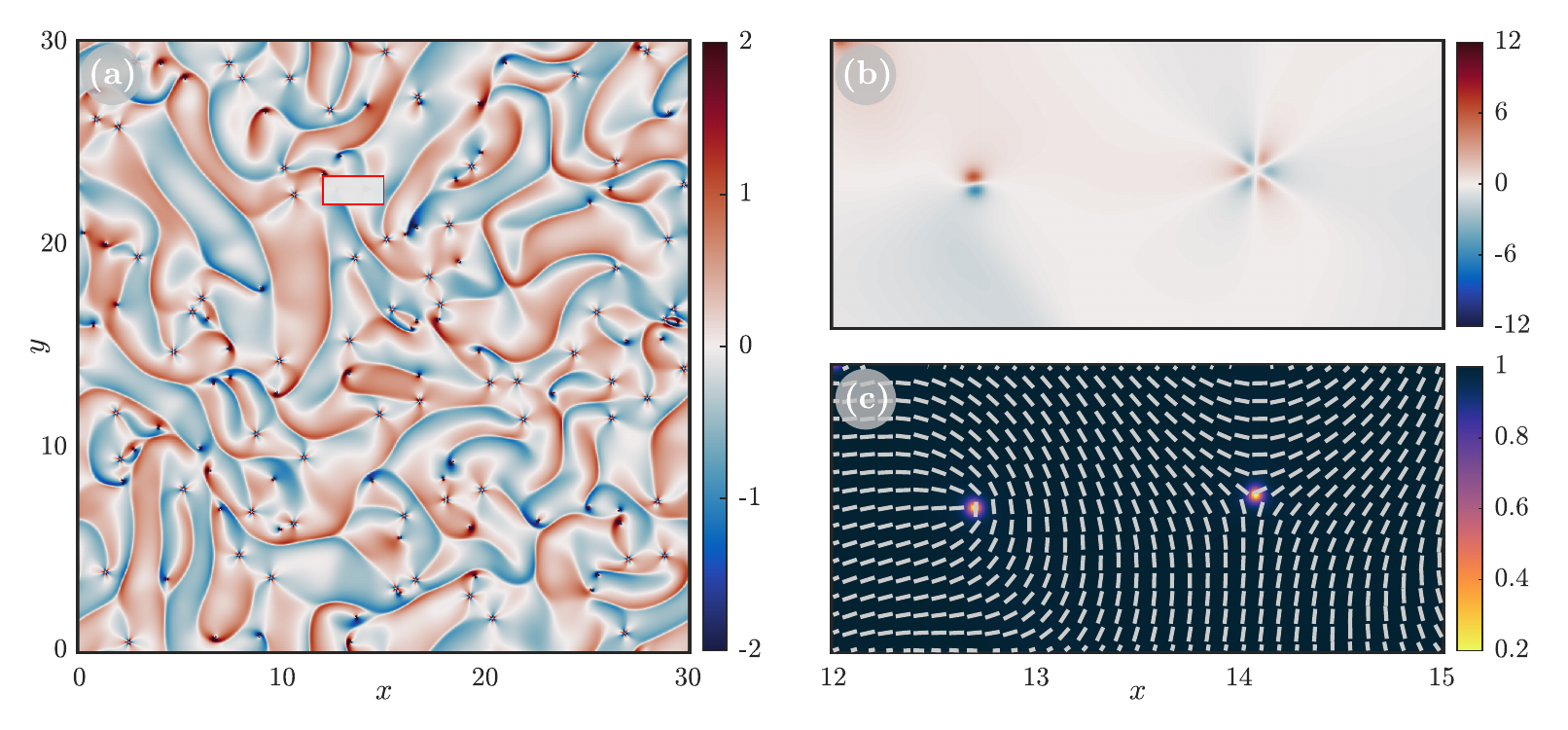}
  \vspace{-0.25in}
  \caption{Snapshot of a two-dimensional simulation with alignment strength $\zeta = 64$ and box size $L = 30$. Panel (a) shows the vorticity field, which consists of isolated vortex dipoles and hexapoles and large regions of nearly constant vorticity. A close up of the vorticity and scalar order fields, indicated by the small box outlined in red in panel (a), shows a matching of symmetry between the dipoles/hexapoles and the $\pm1/2$ disclinations. (Movies of this and the following simulation(s) can be found in the supplementary material.)}\label{fig:steric-2d}
  \end{figure}

    Figure \ref{fig:steric-2d} shows a snapshot of a two-dimensional simulation at a late time for $\zeta = 64$ where the time step is $\Delta t = 10^{-4}$. The vorticity field, shown in panel (a), consists of isolated vortices which trail shock-like structures in the global field. Close ups of the scalar order and vorticity fields near two defects are shown in panels (b) and (c), which show dipole and hexapole structures whose 1- and 3-fold symmetries are inherited by the $\pm1/2$ sign of the defect, respectively.

\begin{figure}[t!]
  \centering
  \includegraphics{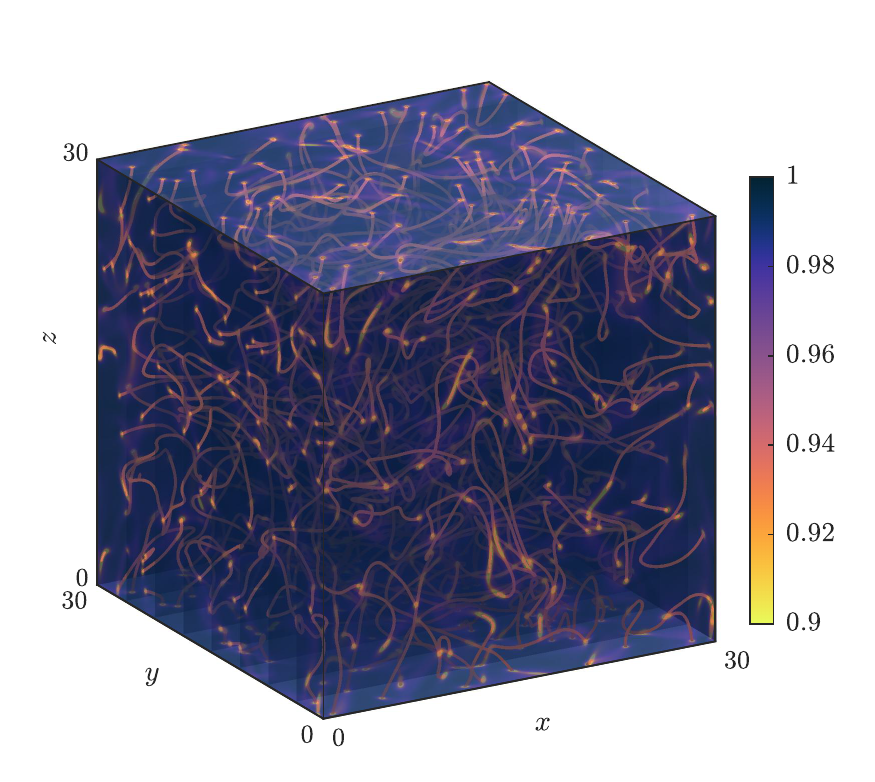}
  \includegraphics{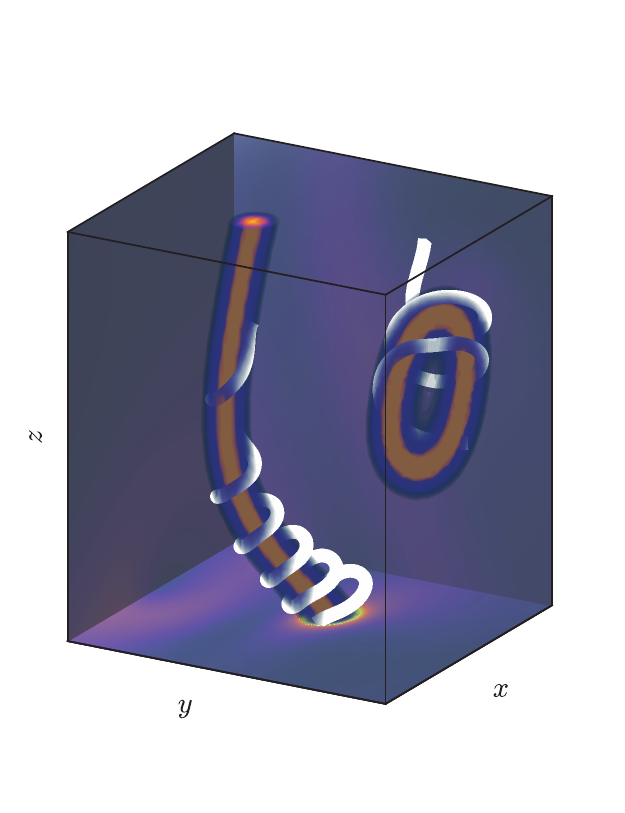}
  \vspace{-0.125in}
  \caption{Three-dimensional simulation with alignment strength $\zeta = 8$ and box size $L = 30$. Panel (a) shows three-dimensional contours of the scalar order field, which reveal relatively isolated disclination lines and loops. A close up of these disclinations is shown in panel (b), with example vortex field lines (white) superimposed. The vortex lines wind around the disclinations, similar to mutual interactions between vortex lines in Navier-Stokes turbulence.}\label{fig:steric-3d}
  \end{figure}

    A three-dimensional simulation is shown in Figure \ref{fig:steric-3d} for $\zeta = 8$, where the time step is $\Delta t = 0.0025$. We find the scalar order field, shown in panel (a), consists primarily of long tubes of low orientational order, which is consistent with simulations of the Landau-deGennes theory \cite{Duclos:2020}. The three-dimensional vortex field lines, shown near isolated disclination lines in panel (b), wind around the axis of the disclination. Such intertwining structures are observed in vorticity in classical three-dimensional turbulence, and are the analogous extension of the dipoles observed in the previous two-dimensional simulation.

    \begin{figure}[t!]
    \centering
    \includegraphics{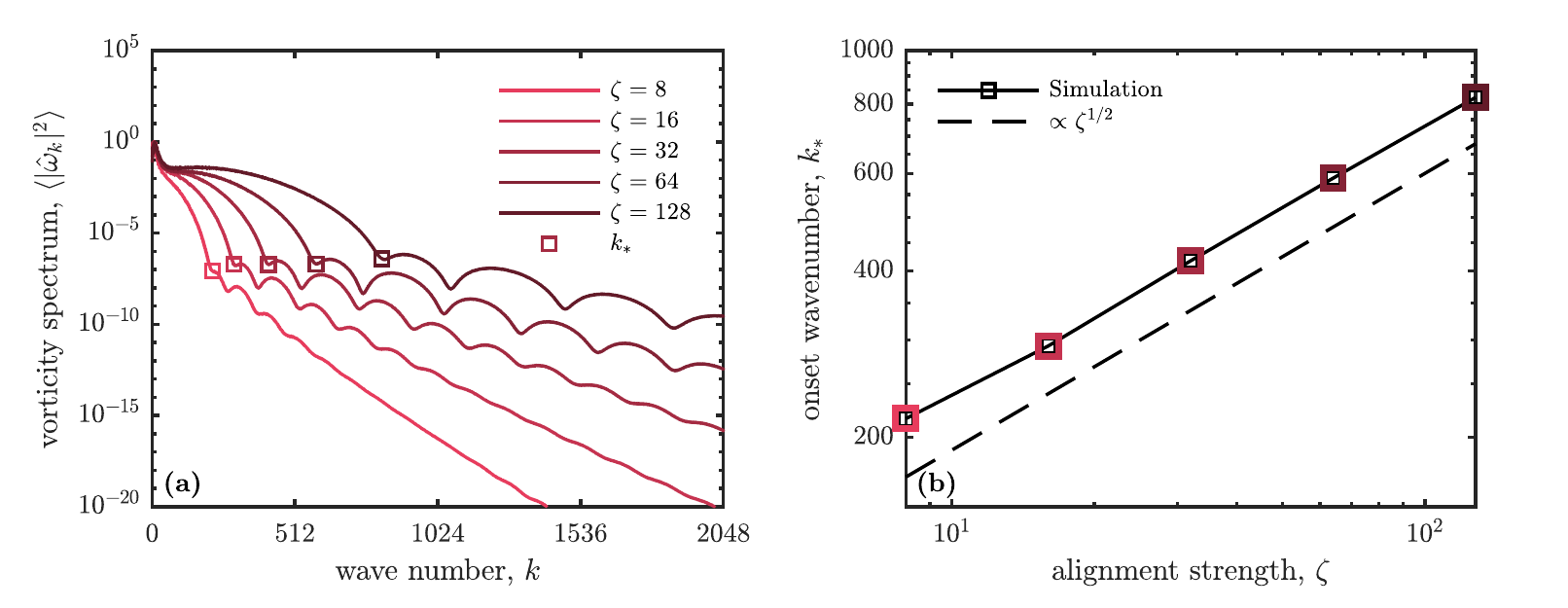}
    \vspace{-0.25in}
    \caption{Characterization of the vortex length scale at strong alignment. Panel (a) shows the squared magnitude of the vorticity spectra for several values of $\zeta$. The oscillations reveal a length scale which can be associated with the radius of decay of the vortex dipoles and hexapoles generated by disclinations in the director field. Panel (b) shows the onset wavenumber $k_*$ of the oscillations (identified by the squares in panel (a)) against the alignment strength, in which we find a $k_*\propto \zeta^{1/2}$ scaling, which is consistent with analytical predictions.}\label{fig:steric-scaling}
    \end{figure}
    
    These and further simulations show that the vortex structures above exhibit length scales that are strongly coupled to the alignment parameter $\zeta$. We can get an analytical estimate on these length scales by rescaling the coarse-grained equations (\ref{eq:stress-nd})-(\ref{eq:dD/dt}). Defining $t' = \zeta t$, $\x' = \zeta^{1/2}\x$ and $\u' = \zeta^{-1/2}\u$, we find
    \begin{gather*}
    -\Delta'\u' + \grad'q' = \grad'\cdot\bm\Sigma',\\
    \grad'\cdot\u' = 0,
    \end{gather*}
    and
    \begin{gather*}
    \D^{\grad'} + 2\S:\E' = 4(\D\cdot\D - \S:\D) + d_T\Delta'\D - 2dd_R\zeta^{-1}\bpar{\D - \frac{c}{d}\I},
    \end{gather*} 
    where the rescaled stress is
    \[ \bm\Sigma' = \alpha\zeta^{-1}\D + \beta\S:\E' - 2\beta(\D\cdot\D - \S:\D).\]
    In the limit $2d_R/\zeta \ll 1$ and $\alpha/\zeta \ll 1$, this system of equations becomes independent of $\zeta$ so that the characteristic length $\ell_c$ must scale as $\ell_c \propto \zeta^{-1/2}$, regardless of the spatial dimension.

    To assess this asymptotic regime, we run several two-dimensional simulations with successively doubled values of $\zeta$ and compare the vorticity spectra at a late time. As shown in Figure \ref{fig:steric-scaling}a, the instantaneous spectra exhibit regular oscillations that increase in amplitude and width as $\zeta$ increases. Notably, the spectra reveal a length scale corresponding to the onset wave number $k_*$ at which the oscillations begin. Figure \ref{fig:steric-scaling}b shows this onset wavenumber as a function of $\zeta$, which indicates a $k_*\propto\zeta^{1/2}$ scaling, in agreement with the analytical prediction as the wave number has units of inverse length. Taken with the predicted characteristic time scale $t_c\propto \zeta^{-1}$, this scaling could also be used to characterize the number and rate of creation of defects, as well as their typical velocities \cite{Thampi:2014}.


    \subsection{Turbulent dynamics}
    
    A peculiar property of the continuum kinetic model is that linear stability analysis in periodic geometries predicts the smallest wavenumbers are the most unstable \cite{ESS:2013}. Because of this, the linear theory does not predict a characteristic length scale. Nonlinear simulations, however, can provide insight into characteristic length scales in the system and the transfer of energy across them. Here we simulate the coarse-grained model with the Bingham closure to study this nonlinear behavior for large box sizes. Large box simulations allow for more unstable low wave numbers in the system, which we expect to drive increasingly turbulent dynamics. The precise statistics of this so-called active turbulence has been the focus of several recent studies with Landau-deGennes type theories  \cite{Alert:2020,Carenza:2020}, and our formulation allows us to study these statistics with a first-principles approach. For the following simulations we fix the dimensionless parameters $\alpha = -1$, $\beta = 0.8$, $\zeta = 1$, and $d_T = d_R = 0.05$, and vary the linear dimension $L$.

\begin{figure}[t!]
    \centering
    \includegraphics{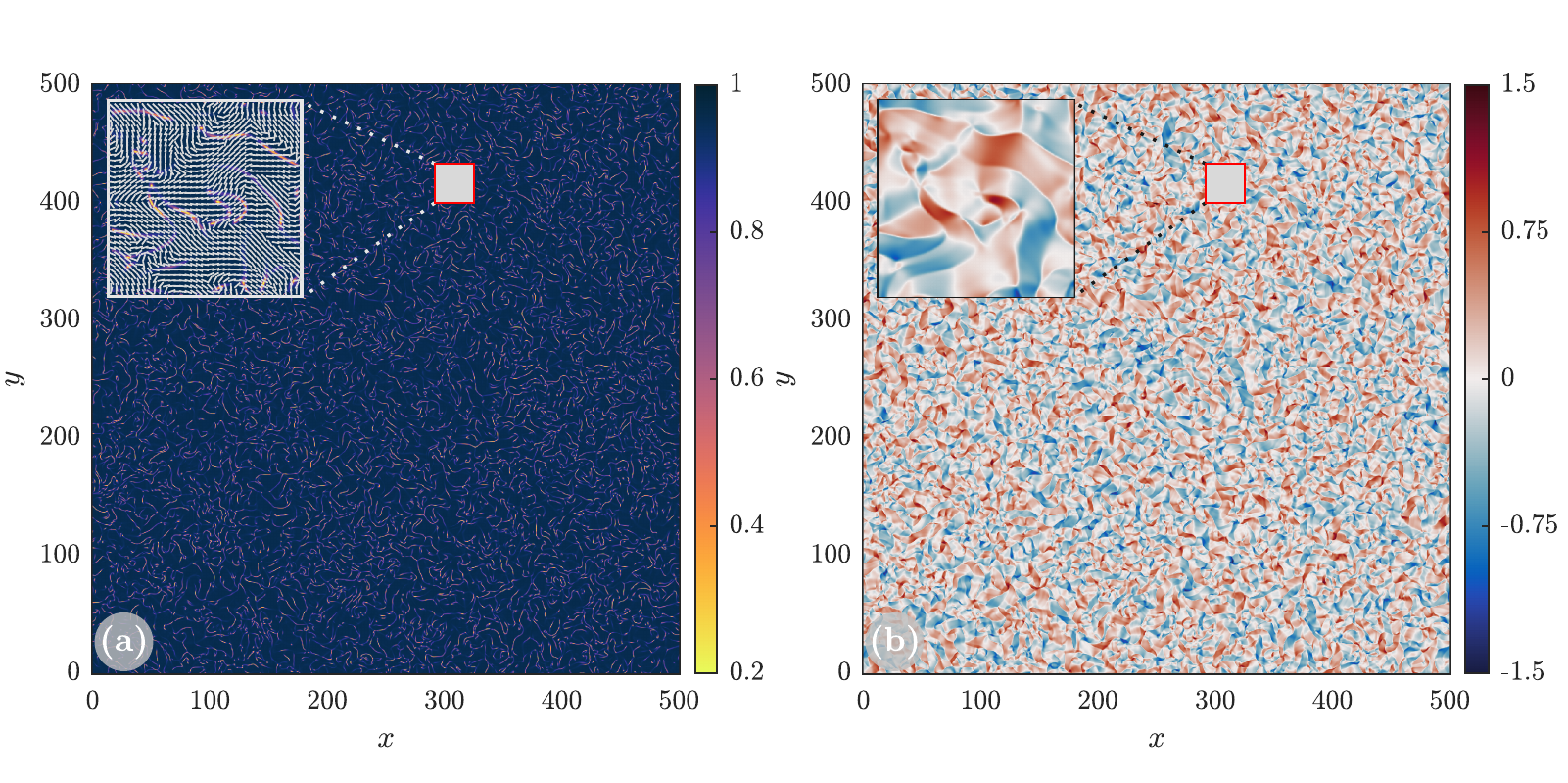}
    \vspace{-0.1in}
    \caption{(a) Scalar order and (b) vorticity fields from a two-dimensional simulation with box size $L = 500$ and alignment strength $\zeta = 1$. The scalar order field consists of many regions of low orientational order connected with topological defects in the director field. In comparison with panel (b), we find the vorticity rapidly changes sign across the bands of low orientational order generating patch-like patterns.}\label{fig:large-box-2d}
    \end{figure}

\begin{figure}[t!]
    \centering
    \includegraphics{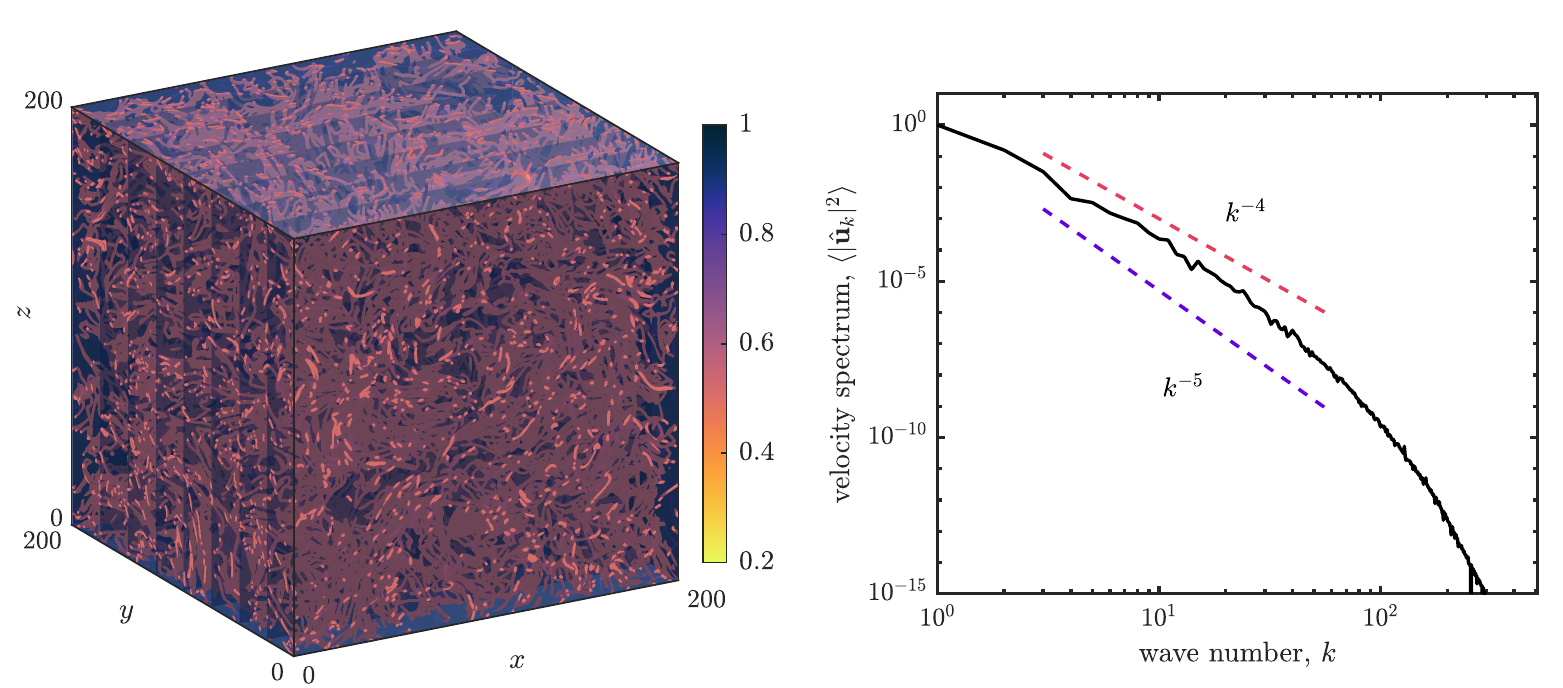}
    \caption{Three-dimensional simulation with box size $L = 200$ and alignment strength $\zeta = 1$. Panel (a) shows contours of the scalar order field, which reveal a dense concentration of disclinations. The motion of these disclinations is coupled to large-scale motion. A rough characterization of this coupling is shown by panel (b), which reveals an approximate power law scaling in the velocity spectrum at low wave numbers, indicating a transfer of energy across length scales in the system.}\label{fig:large-box-3d}
    \end{figure}

    Figure \ref{fig:large-box-2d} shows a snapshot of the scalar order and vorticity fields from a two-dimensional simulation with box size $L = 500$. In contrast to the simulations with strong alignment in Section \ref{sec:strongly-aligned}, the scalar order field is densely packed with topological defects that undergo rapid nucleation and annihilation events. These defects do not seem to create strong vortex dipoles, rather the bands of low orientational order connecting them generate small patches of nearly constant vorticity. In three dimensions, with $L = 200$, we find the scalar order field also consists of fine-scale defect structures, shown in Figure \ref{fig:large-box-3d}a, with many intertwining disclination loops and tubes which also undergo frequent nucleation and annihilation events. 

    In both of these simulations the dense defect structures drive large scale motion. This transfer across scales is often observed in turbulent fluids and can be characterized by analyzing the squared velocity spectrum, which in classical turbulence reflects kinetic energy. (Note that due to the low Reynolds number the kinetic energy of our system has no relevance. However, based on the entropy identity (\ref{eq:dS/dt}), velocity gradients characterize entropy production or dissipation as they would in classical turbulence.) Panel (b) in Figure \ref{fig:large-box-3d} shows the computed velocity spectrum $\langle |\hat\u_k|^2 \rangle$ summed over spherical shells in the wave number $|\k| = k$. Unlike the strongly aligned case the spectrum does not exhibit oscillations, rather, at lower wave numbers we observe an approximate power law between $k^{-4}$ and $k^{-5}$, which transitions to a more rapid decay at $k \approx 80$. This transition wave number may indicate a characteristic turbulent length scale, whose precise interpretation is the subject of future investigation.

\section{Concluding remarks}

We developed a robust numerical method for simulating coarse-grained models of apolar particle suspensions with the Bingham closure. Unlike previous approaches, the closure map is constructed by solving for the Bingham distribution over the entire feasible domain of the second moment tensor $\tilde\D$. By transforming this domain to a square domain, we were able to represent the mapping $\tilde\D\mapsto\tilde \S_B$ by a Chebyshev interpolant for efficient use in simulations. This Chebyshev representation reconstructs the closure to near machine precision, with accuracy that can be finely controlled and balanced against cost by modifying the degree of the interpolant. We found accuracy of the closure map is essential for maintaining spatial convergence in the underlying discretization and resolving high wave number behavior, which was shown to be a consequence of the rotation-based approach. The simulations in Section \ref{sec:convergence} underscore the importance of this fact. In particular, defects in the director field correspond to points where the eigendecomposition of the second moment tensor is ill-defined and the rotation-based approach fails. Such states are fundamental to the underlying physics, and accurately resolving them is essential for retaining the overall structure and statistics of the physical system.

Using this method, we studied regimes of strong alignment and large system size. When alignment is strong, we found coupling between defects in the director field and fluid vorticity, which were consistent with the Landau-deGennes theory. This connection could be used to construct reduced models of defect systems, possibly describing defects in an analogous way to interacting point vortices or vortex filaments in the incompressible Euler equations \cite{Vortex-Methods}. We also analytically derived a scaling law $\ell_c\propto\zeta^{-1/2}$ for the defect length scale, which was confirmed through high resolution two-dimensional simulations. In contrast, for large system size we found the dynamics were turbulent, exhibiting chaotic motion from the defect to system scales. Analyzing the velocity spectrum here revealed an approximate power law scaling at low wave numbers, which may reflect a transfer of energy across length scales in the system. Future work could characterize this transfer of energy more precisely, including its dependence on the system size, the nematic alignment strength, and the magnitude and sign of the active stress.

As formulated here, the Bingham closure only applies to apolar suspensions. In reality, many physical systems are inherently polar, such as microtubule and motor protein assemblies or collections of motile bacteria \cite{GBGBS:2015,Dombrowski:2004}. The Bingham distribution can be generalized to account for polarity, and we are working on similar methods to those developed here to accurately and efficiently construct the generalized closure map.

A significant property of the Bingham closure, which we proved in the Appendix, is that it preserves the evolution of the system entropy, where the entropy is approximated in terms of the Bingham distribution. Combined with the accuracy and efficiency of the method presented here, the Bingham closure could be used to study energetic properties of active systems, particularly in three dimensions, that are consistent with the kinetic theory. 

\section*{Declaration of competing interests}

The authors declare that they have no known competing financial interests or personal relationships that could have appeared to influence the work reported in this paper.

\section*{Author contributions}

{\bf SW}: Formal analysis, investigation, validation, visualization. Writing - original draft, review \& editing. {\bf DBS}: Formal analysis, investigation, validation, visualization. Writing - review \& editing. {\bf MJS}: Formal analysis, investigation, validation, visualization. Writing - review \& editing.

\section*{Acknowledgements}

During the review of this paper, we became aware of work that uses similar methods with a focus towards the Doi theory of passive rod suspensions \cite{Jiang:2021}. We thank Sebastian F\"urthauer for useful discussions. SW acknowledges support from the NSF-GRFP under Grant No. 1839302. MJS acknowledges support by the National Science Foundation under awards DMR-
1420073 (NYU MRSEC) and DMR-2004469.

\bibliography{closure}

\appendix

\section{Calculation of the Bingham parameters}

In the interpolation approach the Bingham distribution is never actually constructed, however its parameters may be needed to compute higher order moments or system statistics, such as the conformational entropy. Here we show how to construct the Bingham parameter $\B$ analytically from the second and fourth moments $\D$ and $\S_B$. First, we compute
\begin{align*}
\int_{|\p|=1} \p\grad_p\Psi_B ~ d\p &=
\int_{|\p|=1} (\p\I - \p\p\p)\cdot\partial_{\p}\Psi_B ~ d\p \\
&= 2\int_{|\p|=1} [(\p\I - \p\p\p)\cdot(\B\cdot\p)]\Psi_B~ d\p \\
&= 2\int_{|\p|=1} (\p\p\cdot\B - \p\p\p\p:\B)\Psi_B~ d\p \\
&= 2(\D\cdot\B - \S_B:\B).
\end{align*}
Integrating by parts gives
\begin{align*}
\int_{|\p|=1} \p\grad_p\Psi_B ~ d\p
&= -\int_{|\p|=1} \grad_p\p ~ \Psi_B ~ d\p \\
&= -\int_{|\p|=1} (\I - d\p\p)\Psi_B ~ d\p \\
&= c\I - d\D,
\end{align*}
which implies
\begin{equation} \D\cdot\B - \S_B:\B = \frac{d}{2}\bpar{\D - (c/d)\I}.\label{eq:ibp}\end{equation}
With $\D$ and $\S_B$ known, this system can be inverted for $\B$. Higher order moments can similarly be determined analytically by integrating by parts with higher order products of $\p$.

\DeclarePairedDelimiter{\moment}{\langle}{\rangle}

\section{Entropy production}

Here we show the Bingham closure satisfies the same energy identity as the kinetic theory \cite{GBJS:2017}, with the entropy represented in terms of the Bingham distribution. For simplicity we assume the concentration is uniform $c(\x,t) = 1$. Throughout we denote $\grad$ as the spatial gradient and $\grad_p = (\I - \p\p)\cdot\partial_\p$ as the gradient operator on the unit sphere. All spatial integrals are assumed to be over the volume $V$.

The steric contribution $\mathcal D(t) = \int (\D - \I/d):(\D - \I/d) ~ d\x$ is only represented in coarse-grained variables which, based on the evolution equation (\ref{eq:dD/dt}), automatically satisfies the same equation for $\mathcal D'(t)$ in both the Bingham closure and the kinetic theory. After some standard manipulations we can show
\begin{align*}
\mathcal D'(t) &= -4dd_R \int (\D-\I/d):(\D-\I/d) ~ d\x + 8\zeta \int \D:(\D\cdot\D - \S_B:\D) ~ d\x \\ & \quad + 4\int \E:(\D\cdot\D - \S_B:\D) ~ d\x - 2d_T \int |\grad\D|^2 ~ d\x.
\end{align*}
Now let $\Psi_B(\x,\p,t) = e^{\gamma(\x,t) + \B(\x,t):\p\p}$ be the Bingham distribution, where $\gamma(\x,t) = -\log Z(\x,t)$ is a normalization factor enforcing $\int_{|\p|=1} \Psi_B d\p = 1$. In terms of $\Psi_B$ the conformational entropy $\mathcal S(t) = \int\int_{|\p|=1} (\Psi/\Psi_0)\log(\Psi/\Psi_0) ~ d\p d\x$ is
\[ \mathcal S(t) = \frac{1}{\Psi_0}\int (\gamma - \gamma_0) + \B:\D ~ d\x,\]
where $\gamma_0 = \log\Psi_0$. Differentiating the constraint $\int \Psi_B ~ d\p d\x = V$ in time gives $\int \gamma_t + \B_t:\D ~d\x = 0$, which implies
\[ \mathcal S'(t) = \frac{1}{\Psi_0}\int \B:\D_t ~ d\x.\]
Using Equation (\ref{eq:dD/dt}) for $\D_t$ we get
\begin{align*}
\B:\D_t &= -\B:\u\cdot\grad\D + \B:(\grad\u\cdot\D + \D\cdot\grad\u^T) - 2\B:(\S_B:\E) 
\\ & \quad +
4\zeta \B:(\D\cdot\D - \S_B:\D) + d_T\B:\Delta\D - 2dd_R(\D - \I/d).
\end{align*}
Contracting the integration by parts identity (\ref{eq:ibp}) against $\E$ and $\D$, respectively, gives
\begin{align*}
\B:(\grad\u\cdot\D + \D\cdot\grad\u^T) - 2\B:(\S_B:\E) = d\D:\E
\end{align*}
and
\begin{align*} 
\B:(\D\cdot\D - \S_B:\D) &= \frac{d}{2}\bpar{\D:(\D-\I/d)}
\\ & =
\frac{d}{2}\bpar{(\D - \I/d):(\D - \I/d)},
\end{align*}
so that
\[ \B:\D_t = -\B:\u\cdot\grad\D + d\D:\E + 2d\zeta (\D-\I/d):(\D - \I/d) + d_T\B:\Delta\D - 2dd_R(\D - \I/d).\]
From the condition $\int \grad\cdot(\u c) ~d\x = 0$ we have $\int\u\cdot\grad\gamma ~ d\x = -\int \u\cdot(\grad\B:\D) ~ d\x $ which, after a few integrations by parts, gives
\[\int\B:(\u\cdot\grad\D) ~ d\x = 0.\]
So far the evolution of the conformational entropy is
\[ \Psi_0 \mathcal S'(t) = d\int \D:\E~ d\x + 2d\zeta\int (\D-\I/d):(\D - \I/d) ~ d\x + d_T\int\B:\Delta\D~ d\x - 2dd_R\int\B:(\D - \I/d) ~ d\x,\]
which we want to write in terms of definitely signed quantities. Multiplying the Stokes equation (\ref{eq:stokes-nd}) by $\E$ and integrating by parts gives
\[ 2\int \E:\E ~ d\x = -\int \E:\bm\Sigma_B ~ d\x,\]
which implies
\begin{align*}
2\int \E:\E ~ d\x &=
-\alpha\int\D:\E ~ d\x - \beta\int\E:\S_B:\E ~ d\x + 2\zeta\beta\int \E:(\D\cdot\D - \S_B:\D) ~ d\x.
\end{align*}
We can use this to solve for $\int\D:\E ~ d\x$,
\begin{align*}
\Psi_0 \mathcal S'(t) &= \frac{d}{\alpha}\bpar{-2\int\E:\E~d\x -\beta\int\E:\S_B:\E ~ d\x + 2\zeta\beta\int\E:(\D\cdot\D - \S_B:\D) ~ d\x} \\ & \quad + 2d\zeta\int (\D - \I/d):(\D - \I/d) ~ d\x + d_T\int\B:\Delta\D~ d\x - 2dd_R\int\B:(\D - \I/d) ~ d\x.
\end{align*}

It is then left to show
\[ \int\B:\Delta\D~ d\x = -\int \int_{|\p|=1} |\grad\log\Psi_B|^2\Psi_B ~ d\p d\x \]
and
\[ \int\B:(\D - \I/d) ~ d\x = \frac{1}{2d}\int\int_{|\p|=1}|\grad_p\log\Psi_B|^2 \Psi_B ~ d\p d\x.\]
For the first term, differentiating in space gives $|\grad\log\Psi_B|^2\Psi_B = |\grad\gamma + \grad\B:\p\p|^2\Psi_B$. Using the condition $\int\grad\cdot(\grad c) ~ d\x= 0$, we find
\begin{align*}
0 &= \int_{|\p|=1} \grad\cdot\grad\Psi_B ~ d\p \\
&= \int_{|\p|=1}\grad\cdot[(\grad\gamma + \grad\B:\p\p)\Psi_B] ~ d\p\\
&= \int_{|\p|=1} (\Delta\gamma + \Delta\B:\p\p)\Psi_B + |\grad\gamma + \grad\B:\p\p|^2 \Psi_B ~ d\p,
\end{align*}
so that $\int\Delta\gamma + \Delta\B:\D  ~ d\x= -\int\moment{|\grad\gamma + \grad\B:\p\p|^2}~d\x$, which, after two integrations by parts, gives $\int\B:\Delta\D~d\x = -\int|\grad\log\Psi_B|^2\Psi_B ~ d\p d\x$ as desired. For the second term, we have $|\grad_p\log\Psi_B|^2 = |2(\I - \p\p)\cdot\B\cdot\p|^2$ so that
\begin{align*}
\int\int_{|\p|=1} |\grad_p\log\Psi_B|^2\Psi_B ~ d\p d\x
&= 4\int \B:(\B\cdot\D) - \B:(\S_B:\B) ~ d\x\\
&= 2d\int \B:(\D - \I/d) ~ d\x,
\end{align*}
where we used the same integration by parts identity (\ref{eq:ibp}). Finally, incorporating the expression for $\mathcal D'(t)$, we have
\begin{equation*}
\begin{aligned}
\mathcal E'(t) &= -\frac{d}{\alpha\Psi_0}\bpar{\int 2\E:\E  + \beta \E:\S_B:\E ~ d\x}+ (2d\zeta/\Psi_0 - 4dd_R\kappa)\int (\D-\I/d):(\D-\I/d) ~ d\x \\
&\quad + 8\kappa\zeta\int\D:(\D\cdot\D - \S_B:\D)~d\x   -2\kappa d_T \int |\grad\D|^2 ~ d\x \\ & \quad -\bbrack{d_T\int\int_{|\p|=1}\Psi_B|\grad\log\Psi_B|^2 ~ d\p d\x + d_R\int\int_{|\p|=1}\Psi_B|\grad_p\log\Psi_B|^2 ~d\p d\x},
\end{aligned}
\end{equation*}
which is the desired expression. (Note that $\kappa = -d\zeta\beta/2\Psi_0\alpha > 0$ is chosen so that the contribution from $\int \E:(\D\cdot\D - \S_B:\D) ~ d\x$ vanishes.)

\section{Chebyshev coefficients $\mu_1\mapsto \tilde S_{1111}$}

\begin{center}
\resizebox{\columnwidth}{!}{
\begin{tabular}{|c|c|c|c|c|c|c|c|c|c|c|} \hline
 & \large{0} & \large 1 & \large 2 & \large 3 & \large 4 & \large 5 & \large 6 & \large 7 & \large 8 & \large 9 \\ \hline
\large{0} & 0.662433067815903 & 0.305096697570660 & 0.022661293663811 & 0.006929073516477 & 0.002508696495226 & 0.000651686393991 & 0.000003149067870 & -0.000146501968392 & -0.000116594835353 & -0.000052333302979 \\
\large{10} & -0.000005876267519 & 0.000014375002152 & 0.000016332739170 & 0.000010270338585 & 0.000003367708799 & -0.000001122818696 & -0.000002764227720 & -0.000002463455278 & -0.000001369444354 & -0.000000311420682 \\
\large{20} & 0.000000333432499 & 0.000000536642052 & 0.000000445442150 & 0.000000238673323 & 0.000000046766023 & -0.000000069897867 & -0.000000107035811 & -0.000000090051728 & -0.000000050610698 & -0.000000012725701 \\
\large{30} & 0.000000011666027 & 0.000000020862528 & 0.000000019015530 & 0.000000011929798 & 0.000000004355589 & -0.000000001058379 & -0.000000003616565 & -0.000000003863898 & -0.000000002815955 & -0.000000001409792 \\
\large{40} & -0.000000000245577 & 0.000000000439389 & 0.000000000669175 & 0.000000000594818 & 0.000000000385800 & 0.000000000169018 & 0.000000000011082 & -0.000000000071943 & -0.000000000094069 & -0.000000000079964  \\
\large{50} & -0.000000000052121 & -0.000000000025272 & -0.000000000006080 & 0.000000000004532 & 0.000000000008540 & 0.000000000008620 & 0.000000000006960 & 0.000000000004899 & 0.000000000003050 & 0.000000000001588  \\
\large{60} & 0.000000000000502 & -0.000000000000256 & -0.000000000000728 & -0.000000000000942 & -0.000000000000940 & -0.000000000000779 & -0.000000000000531 & -0.000000000000267 & -0.000000000000040 & 0.000000000000116  \\
\large{70} & 0.000000000000195 & 0.000000000000207 & 0.000000000000173 & 0.000000000000118 & 0.000000000000060 & 0.000000000000013 & -0.000000000000019 & -0.000000000000035 & -0.000000000000038 & -0.000000000000032  \\
\large{80} & -0.000000000000023 & -0.000000000000013 & -0.000000000000005 & 0.000000000000001 & 0.000000000000004 & 0.000000000000006 & 0.000000000000006 & 0.000000000000004 & 0.000000000000003 & 0.000000000000002  \\
\large{90} & 0.000000000000001 & -0.000000000000000 & -0.000000000000001 & -0.000000000000001 & -0.000000000000001 & -0.000000000000001 & -0.000000000000001 & -0.000000000000000 & -0.000000000000000 & -0.000000000000000 \\ \hline
\end{tabular}
}
\end{center}

\end{document}